\pdfoutput=1

\documentclass[alpha-refs]{wiley-article}

\usepackage{graphicx}
\usepackage[space]{grffile}
\usepackage{latexsym}
\usepackage{textcomp}
\usepackage{longtable}
\usepackage{tabulary}
\usepackage{booktabs,array,multirow}
\usepackage{amsfonts,amsmath,amssymb}
\usepackage{natbib}
\usepackage{url}
\usepackage{hyperref}
\hypersetup{colorlinks=false,pdfborder={0 0 0}}
\usepackage{etoolbox}
\newif\iflatexml\latexmlfalse

\AtBeginDocument{\DeclareGraphicsExtensions{.pdf,.PDF,.eps,.EPS,.png,.PNG,.tif,.TIF,.jpg,.JPG,.jpeg,.JPEG}}

\usepackage[utf8]{inputenc}
\usepackage[english]{babel}

\usepackage{siunitx}

\usepackage[T1]{fontenc} 

\iflatexml

\else

\paperfield{causal inference}
\corraddress{Jaffer M. Zaidi, Department of Biostatistics, Univerity of North Carolina, Chapel Hill, 27599, USA}
\corremail{jzaidi@email.unc.edu}
\fundinginfo{National Institutes of Health,\\
Grant/Award Number: 5T32AI007358-30;\\
National Science Foundation\\
DMS-1127914}
\fi

\papertype{Original Article}

\title{On the identification of individual level pleiotropic, pure direct, and
principal stratum direct effects without cross world assumptions}

\author[1]{Jaffer M. Zaidi}

\affil[1]{Department of Biostatistics, University of North Carolina}

\author[2]{Tyler J. VanderWeele}

\affil[2]{Department of Biostatistics and Epidemiology, Harvard T.H. Chan School of Public Health}

\runningauthor{Jaffer Zaidi}

\begin{document}

\maketitle
\selectlanguage{english}
\begin{abstract}
The analysis of natural direct and principal stratum direct effects has
a controversial history in statistics and causal inference as these
effects are commonly identified with either untestable cross world
independence or graphical assumptions. This paper demonstrates that the
presence of individual level natural direct and principal stratum direct
effects can be identified without cross world independence assumptions.
We also define a new type of causal effect, called pleiotropy, that is
of interest in genomics, and provide empirical conditions to detect such
an effect as well. Our results are applicable for all types of
distributions concerning the mediator and outcome.%
\end{abstract}%

\textbf{Keywords} --- direct effects, mediation, missing data, natural direct
effects, pleiotropy, principal stratum.

\section{Introduction}

{\label{707961}}

\subsection{Natural and Principal Stratum Direct
Effects}

{\label{485113}}

~~~~The analysis of natural direct and principal stratum direct effects
has a long history in statistics and causal
inference~\citep{pearl2009causality,vanderweele2015explanation,frangakis2002principal,robins1992identifiability}. All of these effects are defined in a
counterfactual or potential outcome framework ~\citep{pearl2009causality,vanderweele2015explanation,frangakis2002principal,robins1992identifiability}, but
these effects are commonly identified using either untestable cross
world independence assumptions or graphical
assumptions~\citep{pearl2009causality,vanderweele2015explanation,frangakis2002principal,robins1992identifiability}. Consequently, the identification and use
of these effects in research has elicited controversy from authorities
across varied disciplines~\citep{pearl2009causality,vanderweele2015explanation,frangakis2002principal,robins1992identifiability}. The invited discussion
section of a recent paper elucidates some of the controversy surrounding
natural direct and indirect effects~\citep{imai2013experimental}.~

~~~~We detail a few of the issues surrounding natural direct effects.
The first concern is that these effects are identified under
unverifiable and untenable assumptions that usually offer no
experimental design to detect and quantify individual level direct and
indirect effects. The second grievance is that the cross world
counterfactuals that are used to define natural direct and indirect
effects might not even exist for all or a subset of the population, and
therefore using cross world counterfactuals that might not even have a metaphysical
interpretation or even exist in the simplest thought experiments
presents a deep challenge to scientists and statisticians who are
concerned with interpreting the results of actual experiments.
Statisticians have attempted to overcome these first two challenges
through different approaches~\citep{imai2013experimental}. In the online supplement, the first author details a few of his concerns regarding natural direct and indirect effects that are not shared with the second author.

~~~~The first contribution of this paper demonstrates that individual
level natural direct and principal stratum direct effects can sometimes
be empirically detected and quantified without cross world independence
assumptions. We first provide the simpler theory and methodology to
quantify individual level direct effects in the setting of binary
variables, and then follow up with the more complex theory and methods
to quantify unit level direct effects in the setting of tuple or vector
valued continuous mediator and/or outcome. The second contribution of
this paper is deriving lower bounds on the magnitude of the detected
direct effect in comparison to the other causal effects. Such bounds~
enable statisticians to uncover mechanisms of the effect of a treatment
on an outcome. Related results when one assumes that the treatment has a
monotonic effect on the potential mediator are also derived. We also
offer causally interpretable and user friendly sensitivity analysis for
assessing the consequences of violations of the monotonicity
assumption.~ 

~~~~As a consequence of our work, the quantification of individual level
direct effects can sometimes be carried out with no stronger
identification assumptions than the assumptions needed to analyze total
effects in randomized clinical trials~\citep{vanderweele2015explanation,imbens2015causal}. Finally, we
also define ``individual level pleiotropic effects'' to be present when
a treatment is shown to cause two different outcomes in a single
individual in the population. Such an effect is of interest in
understanding the etiology of different phenotypes from single gene or
allele.

\par\null

\section{Quantification of individual level direct
effects}

{\label{815563}}

\subsection{Notation and assumptions}

{\label{481406}}

~~~~Let~\(Y\) denote a binary outcome. In addition, suppose
a binary treatment~\(X\in\{0,1\}\) is randomized at baseline, and a
binary variable~\(M\in\{0,1\}\) is also observed at a time between
measurement of ~\(Y\)~and the assignment of
treatment~\(X\). The individuals or units in our study,
denoted by symbol~\(\omega\) compose a finite population,
denoted~\(\Omega.\)~Define counterfactuals~\(Y_{x}(\omega)\)
and~\(M_{x}(\omega)\) to be the value
of~\(Y\) and~\(M\) respectively had we set the
value of the treatment~\(X\) to~\(x\) for
individual~\(\omega\). The manipulated
counterfactual~\(Y_{x,m}(\omega)\) denotes the value
of~\(Y(\omega)\) had we set~\(X\)
to~\(x\) and~\(M\) to~\(m\).
Finally, define the counterfactual~\(Y_{x,M_{x^{\prime}}}(\omega)\) to be the value
of~\(Y\) had we set the value of the
treatment~\(X\) to~\(x\) and set the value
of the variable~\(M\) to the value it would take had we
set~\(X\) to some value~\(x^{\prime}\), that could be
different to~\(x\) for individual~\(\omega\).
Denote~\(P_c(y_1,y_0,m_1,m_0)=P(Y_1=y_1,Y_0=y_0,M_1=m_1,M_0=m_0),\) which under the composition axiom is equal to
\(P(Y_{1,M_1}=y_1,Y_{0,M_0}=y_0,M_1=m_1,M_0=m_0)\).

~~~~The individual level total effect is defined
as~\(Y_{1}(\omega)-Y_{0}(\omega)\).~The individual pure direct effect is defined
as~\(Y_{1,M_{0}}(\omega)-Y_{0,M_{0}}(\omega)\)~and the total indirect is defined
as~\(Y_{1,M_{1}}(\omega)-Y_{1,M_{0}}(\omega)\)~\citep{robins1992identifiability}. Similarly, the total direct
is defined as~\(Y_{1,M_{1}}(\omega)-Y_{0,M_{1}}(\omega),\) and the individual pure indirect
effect is defined as~\(Y_{0,M_{1}}(\omega)-Y_{0,M_{0}}(\omega)\)~\citep{robins1992identifiability}. These pure
and total direct and indirect effects are sometimes also referred to as
natural direct and indirect effects~\citep{pearl2009causality}. The individual
controlled direct effect is defined as~\(Y_{1,m}(\omega)-Y_{0,m}(\omega)\). We say there
is an individual principal stratum direct effect if for some
individual~\(\omega\) for whom ~\(M_{1}(\omega)=M_{0}(\omega)=m,\) we
have~\(Y_1(\omega)-Y_0(\omega)\ne0.\) Population total, total direct, pure direct,
total indirect, pure indirect and controlled direct effects are derived
through taking expectations of the relevant individual level effects.
The population level principal stratum direct effect of
stratum~\(M_{1}(\omega)=M_{0}(\omega)=m\) is defined as~\(E(Y_{1}-Y_{0}\mid M_{1}=M_{0}=m)\) for
some~\(m\in\{0,1\}\)~\citep{frangakis2002principal}. For
binary~\(M\) and~\(Y,\) the unstandardized
principal stratum direct effects for some~\(m\in\left\{0,1\right\}\)
and~\(y\in\left\{0,1\right\}\) is given as \(P\left(Y_1=y,Y_0=1-y,M_1=m,M_0=m\right)-P\left(Y_1=1-y,Y_0=y,M_1=m,M_0=m\right).\)The corresponding
standardized principal stratum direct effects for some~\(y\in\left\{0,1\right\}\)
and~\(m\in\left\{0,1\right\}\) is given as~
$$\frac{P\left(Y_1=y,Y_0=1-y,M_1=m,M_0=m\right)-P\left(Y_1=1-y,Y_0=y,M_1=m,M_0=m\right)}{P\left(M_1=m,M_{0\ }=m\right)}.\nonumber$$

~~~~The proportion of individuals for whom treatment causes (or
prevents) the outcome among the principal stratum~\(M_1\left(\omega\right)=M_0\left(\omega\right)=m\) for~\(y=1\) (\(y=0\) respectively) is denoted \(P\left(Y_1=y,Y_0=1-y\mid M_1=m,M_0=m\right).\) Each of the principal stratum direct effects
presented herein differ from the usual definition of~\(E\left(Y_1-Y_0\mid M_1=M_0=m\right)\)
provided by~\citet{frangakis2002principal}. Only in the case of a binary outcome does the usual \citet{frangakis2002principal} definition of a principal stratum direct effect \(E\left(Y_1-Y_0\mid S_1=S_0=1\right)\) and the standardized principal direct effect \(P\left(Y_1=y,Y_0=1-y\mid M_1=m,M_0=m\right)-P\left(Y_1=1-y,Y_0=y\mid M_1=m,M_0=m\right)\) coincide with each other. We emphasize each type of causal
estimand to highlight the relative dissimilarity and relations between
each of these population level causal estimands for researchers to
understand and deliberate on what needs to be estimated. In addition to
the question of which causal estimand is of practical and scientific
relevance, researchers will also need to ask how are these causal
contrasts estimated solely from the observed data distribution?

~~~~With treatment randomized at baseline, one may make the `weak
ignorability' assumption~\((Y_{x},M_{x})\amalg X\). This is the only assumption that is guaranteed by design of a randomized trial. Unlike typical approaches
to identify natural direct and indirect effects, no cross world
independence assumptions of the
form~\(Y_{x,m} \amalg M_{x^{\prime}}\mid X\)~\citep{pearl2009causality} are used to derive our
results. Throughout, we require the consistency assumption for
both~\(Y_{x}\) and~\(M_{x}\)~, which means that
when~\(X(\omega)=x\) then~\(Y_{x}(\omega)=Y(\omega)\)
and~\(M_{x}(\omega)=M(\omega)\). This assumption states that the value of
~\(Y\)and~\(M\) that would have been
observed if~\(X\) had been set to what in fact they were
observed to be is equal respectively to the value of~\(Y\)
and~\(M\) that was observed. Additionally, the consistency
assumption for~\(Y_{x,m}\) means that when~\(X(\omega)=x\)
and~\(M(\omega)=m,\) then~\(Y_{x,m}(\omega)=Y(\omega)\). Finally, we assume the
composition axiom, which states that when~\(M_x(\omega)=m\) for
some~\(m\in\{0,1\}\), then~\(Y_{x,m}(\omega)=Y_x(\omega)\). The randomization,
consistency, and composition assumptions are the only assumptions
necessary to derive our results. The assumptions needed to quantify such
direct effects in observational studies is provided in the online
supplement. Proofs of all results are given in the Appendix.

\subsection{Identification of direct effects under randomization}

{\label{801939}}
\par\null
\begin{theorem}
Suppose $X$ is randomized at baseline. If for some $y\in\{0,1\}$ and $m\in\{0,1\},$ $P(Y=y,M=m\mid X=1)+P(Y=1-y,M=m\mid X=0)>1,$ then there exists a non-empty subpopulation $\Omega_{s}\subseteq\Omega$ such that for every $\omega\in\Omega_{s},$ $Y_{1}(\omega)=y,$ $Y_{0}(\omega)=1-y,$ and $M_{1}(\omega)=M_{0}(\omega)=m.$ 
\end{theorem}

\begin{corollary}
Suppose $X$ is randomized at baseline. If the empirical condition $P(Y=1,M=0\mid X=1)+P(Y=0,M=0\mid X=0)>1$ holds, then there exists a non-empty subset $\Omega_s\subseteq\Omega$ such that for every $\omega \in \Omega_{s}$, $Y_{1,M_1}(\omega)-Y_{0,M_0}(\omega) = Y_{1,M_1}(\omega)-Y_{0,M_1}(\omega)=Y_{1,M_0}(\omega)-Y_{0,M_0}(\omega)=Y_1(\omega)-Y_0(\omega) =1,$ with $M_1(\omega)=M_0(\omega)=0,$ and $Y_{1,M_1}(\omega)-Y_{1,M_{0}}(\omega)=Y_{0,M_1}(\omega)-Y_{0,M_{0}}(\omega)=0,$  i.e.  the total effect, total direct effect, pure direct effect, principal stratum direct effect with $M_{1}(\omega)=M_{0}(\omega)=0$ are all equal to $1$, and the total indirect and pure indirect effects for individual $\omega$ are both zero.
\end{corollary}

~~~~Theorem 1 and Corollary 2 allow for the empirical detection of
individual natural and principal stratum direct effects. A somewhat
related set of results to Theorem 1 are the empirical condition to
detect individual level sufficient cause interaction~\citep{vanderweele2008empirical,vanderweele2012general,ramsahai2013probabilistic}.
Theorem 1 demonstrates that such logic can also be used to detect
individual level natural and principal stratum direct effects. Each of
the four empirical conditions provided in Theorem 1, will correctly
identify when an individual principal stratum direct effect is present
in the population. ~Readers familiar with literature on the
falsification of the binary instrumental variable
model~\citep{balke1997bounds} will recognize that each of the empirical
conditions in Theorem 1 as the complement of one of four `instrumental
inequalities'~\citep{balke1997bounds,Swanson_2018}. The result closest to ours is the
result of~\citet{cai2008bounds} who demonstrate that non-zero population
average controlled direct effects can be identified under the `strong
ignorability' assumption through rejecting the `instrumental
inequality'~\citep{cai2008bounds}, but their result does not make any
statement about individual or population principal stratum or natural
direct effects.

~~~~\citet{richardson2010analysis} also discuss implications of falsifying the
instrumental variable model, but their their results focus on the
controlled direct effect and make no mention of individual principal
stratum direct effects.~\citet{sjolander2009bounds} provides bounds for the
population natural direct effect. These bounds on the population natural
direct effect are different to the~\citet{cai2008bounds} bounds on the
population controlled direct effect. It is simple to demonstrate that if
the~\citet{sjolander2009bounds} lower bound (upper bound) of the natural direct
effect is greater than zero (less than zero) respectively, then one of
the inequalities associated with Theorem 1 holds. However, if we
consider the converse, i.e. we can have an inequality with Theorem 1
can hold, but the~\citet{sjolander2009bounds} bounds for the population natural
direct effect could still include zero. The difference stems
from~\citet{sjolander2009bounds} presenting bounds on the population average
natural direct effects, while we present lower bounds on the proportion
of individuals that display a particular form of a principal stratum
direct effect and consequently a natural direct effect.~

~~~~Other authors derive bounds on the principal stratum direct effects
assuming some form of monotonicity assumption and additionally that the
relevant principal stratum exists~\citep{zhang2003estimation,Hudgens2003}, whereas our
results demonstrates that~when the inequality is satisfied the relevant
principal stratum must in fact exist, as shown below, and provides
information on the proportion of individuals that display such a
principal stratum direct effect in comparison to other counterfactual
responses. In comparison to other approaches, we only require
randomization at baseline, and we bound the proportion of individuals
that display a particular direct effect. To the best of our knowledge,
we are the first to detect and quantify the magnitude of an individual
level direct effect with only the assumptions that are guaranteed by a
randomized experiment.

~~~~Natural direct and indirect effects are identified also if one
assumes an underlying nonparametric structural equation models that
imply the cross world independence assumption~\citep{pearl2009causality}. We do
not make that assumption here. Moreover, under a nonparametric
structural model framework, these natural direct effects fail to be
identified if there exists a post baseline confounder
between~\(M\) and~\(Y\) that is also effected
by treatment~\(X\)~\citep{avin2005identifiability}. However, Theorem 1
will still correctly identify non-zero direct effects even in the
presence of a confounder between~\(M\)
and~\(Y\) that is affected by
treatment~\(X.\)~\citet{imai2013experimental} use a more elaborate
experimental design in which both treatment and mediator are randomized
to identify the presence of indirect effects, but we are not assuming
randomization of the mediator, only treatment. We now give an additional
result on the distribution of counterfactual response types
corresponding to direct effects.

\begin{proposition}
Suppose $X$ is randomized at baseline.
For some $y\in\{0,1\}$ and $m\in \{0,1\}$, the expression\\ $P\left(Y=y, M = m \mid X=1\right) + P\left(Y=1-y,M=m\mid X=0\right)-1$ is equal to
\begin{align}
P_c(y,1-y,m,m)  -\{P_c(1-y,y,m,m)+P_c(1-y,1-y,m,1-m)+P_c(1-y,y,m,1-m)\nonumber\\ +P_c(y,y,1-m,m)+P_c(1-y,y,1-m,m)+P_c(y,1-y,1-m,1-m)+\nonumber\\
 P_c(1-y,1-y,1-m,1-m)+P_c(y,y,1-m,1-m)+P_c(1-y,y,1-m,1-m)\}.\nonumber
\end{align}
\end{proposition}
~~~~The proposition above provides the full consequence of rejecting the
instrumental inequality. The following corollary provides researchers
with the ability to bound the magnitude of the direct effect of interest
with only the assumption of randomization at baseline. We note here that
such results are nowhere present in any text or article, despite the
controversy surrounding direct and indirect effects.

\begin{corollary}
Suppose $X$ is randomized at baseline. For some $m\in\{0,1\}$ and $y\in \{0,1\}$, the expression \\ $P\left (Y=y,M=m\mid X=1\right)+P\left(Y=1-y,M=m\mid X=0\right)-1$ is a lower bound of $P_c(y,1-y,m,m)-P_c(1-y,y,m,m).$

If for some $m\in{0,1}$ and $y\in \{0,1\},$  the inequality $P \left (Y=y,M=m\mid X=1\right )+P(Y=1-y,M=m\mid X=0)>1$ holds, then the expression 
$$\frac{P \left (Y=y,M=m\mid X=1\right )+P(Y=1-y,M=m\mid X=0)-1}{P \left (M=m\mid X=1\right )-P(M=1-m\mid X=0)}\nonumber$$
is a lower bound of 
$$\frac{P_c(y,1-y,m,m)-P_c(1-y,y,m,m)}{P(M_1=M_0=m)},\nonumber$$

which relates to the standardized population principal stratum direct effect. This standardized population principal stratum direct effect is equivalent to
$$P(Y_1=y,Y_0=1-y\mid M_1=m,M_0=m)-P(Y_1=1-y,Y_0=y\mid M_1=m,M_0=m).\nonumber$$
\end{corollary}

{~~~~Theorem 1 demonstrates that individual level principal stratum are
identified solely assuming randomization and consistency of
counterfactuals. Corollary 2 then demonstrates that this individual
level principal stratum direct effect is also a natural direct effect
using the composition axiom. Finally, Proposition 3 provides the full equality relationship between the observed world probability contrasts under investigation and counterfactual probabilities that inform mechanism of treatment. This relationship demonstrates that if
for some~}\(y\in \{0,1\}\) and~\(m\in \{0,1\}\) one rejects a null
hypothesis of the form~\(P(Y=y,M=m\mid X=1)+P(Y=1-y,M=m\mid X=0)\leq 1,\) then the counterfactual
distribution of~\((Y_{1,M_1},Y_{0,M_0},M_1,M_0)\) is constrained. This constraint on
the counterfactual distribution has implications for understanding
mechanisms in the population. It informs us that the proportion of
individuals for whom the treatment has a positive (or negative) total
effect on the outcome and simultaneously the treatment does not change
the mediator value from some fixed value ~\(m\in\{0,1\}\) is greater
than the proportion of individuals for whom the treatment has a negative
(or positive respectively) total effect and simultaneously the treatment
does not change the mediator value from some fixed
value~\(m\in\{0,1\}.\)~

~~~~Corollary 4 tells us the magnitude of the resulting direct effect.
From Corollary 4, researchers can use~\(P(Y=1,M=0\mid X=1)+P(Y=0,M=0\mid X=0)-1,\) which is a
function only of the observed data distribution, as a lower bound on the
risk difference of the proportion of~individuals randomized at baseline
for whom the treatment leaves the variable~\(M\)~fixed at 0
and causes the outcome~\(Y\)~and the proportion of
individuals randomized at baseline for whom the treatment also leaves
the variable fixed at 0 and prevents the
outcome~\(Y\).\(\) In
counterfactual notation, this risk difference is given
as~\(P_c\left(1,0,0,0\right)-P_c\left(0,1,0,0\right).\)~From Corollary 4,~ if~\(P\left(Y=1,M=0\mid X=1\right)+P\left(Y=0,M=0\mid X=0\right)>1,\) then~\(\left\{P\left(Y=1,M=0\mid X=1\right)+P\left(Y=0, M=0\mid X=0\right)-1\right\}/\{P\left(M=0\mid X=1\right)-P\left(M=1\mid X=1\right)\}\) is a lower bound on the risk
difference between the proportion of individuals for who treatment
causes the outcome among the individuals in the principal
strata~\(M_1\left(\omega\right)=M_0\left(\omega\right)=0\) and the proportion of individuals for who
treatment prevents the outcome among the principal
strata~\(M_1\left(\omega\right)=M_0\left(\omega\right)=0,\) which is denoted~\(P\left(Y_1=1,Y_0=0\mid M_1=M_0=0\right)-P\left(Y_1=0,Y_0=1\mid M_1=M_0=0\right).\)

~~~~These theorems, propositions and corollaries provided in this
section cannot be improved upon without additional assumptions. In this
section, we only assumed randomization at baseline, and derived
empirical conditions to detect and quantify the magnitude of individual
level direct effects. Corollary 4 also provides the~ greatest~lower
bound on the proportion of individual level direct effects of interest
when only assuming randomization at baseline in terms of the observed
data distribution. To derive lower bounds that are greater or sharper
than those presented in Corollary 4, monotonicity ~assumptions are
required and are presented next.~~

~~~~Note the bound relating to the standardized population principal
stratum direct effect presented in Corollary 4 requires the
that~\(P\left(Y=y,M=m\mid X=1\right)+P(Y=1-y,M=m\mid X=0)>1.\) The reason for this requirement is that we
first need evidence for an individual level principal stratum direct
effect. Otherwise, the statistician has no evidence with from the data
alone that the population average principal stratum direct effect is
non-zero, because if for some~\(y\in\left\{0,1\right\}\)
and~\(m\in\left\{0,1\right\}\)~\(P\left(Y=y,M=m\mid X=1\right)+P(Y=1-y,M=m\mid X=0)\le1,\) then it is entirely possible
that there is no individual such that~\(M_1\left(\omega\right)=M_0\left(\omega\right)=m\)
and~\(Y_1\left(\omega\right)=y,Y_0\left(\omega\right)=1-y\)~ in the observed data. If~\(P\left(Y=y,M=m\mid X=1\right)+P\left(Y=1-y,M=m\mid X=0\right)-1\) is
less than zero, then take zero as the lower bound of~\(P\left(Y_1=y,Y_0=1-y\mid M_1=m,M_0=m\right),\)
since a probability is always nonnegative.

\par\null

\section{Identification of direct effects with monotonicity
assumptions}

~~~~Now consider the setting of positive monotonicity defined
as~\(M_1(\omega)\ge M_0(\omega)\)~ for all~\(\omega\in\Omega,\) that is, there is
no individual for whom treatment prevents the
mediator~\(M\) from occurring. While such monotonic
effects are never verifiable, they are falsifiable and can sometimes be
justified with subject matter knowledge. Under the assumption of
monotonicity, there is a gain in the capacity to detect direct effects
in comparison to the results from the previous section. Additionally, if
such assumptions are reasonable, the tighter bounds enable detection of
two different principal stratum direct effects in a population. We offer
causally interpretable and user friendly sensitivity analysis in the
setting where the statistician believes that some individuals might
violate the monotonicity assumptions to assess the consequence of how
such assumptions can impact study conclusions.

\par\null
\subsection{\texorpdfstring{{Positive monotonicity of exposure on
mediator}}{Positive monotonicity of exposure on mediator}}

{\label{650797}}

\begin{theorem}
Suppose $X$ is randomized at baseline. In addition, suppose that $X$ has a positive
monotonic effect on $M.$ If for some $m\in\{0,1\}$ and $y\in\{0,1\}$ the inequality $P(Y=1-y,M=m\mid X=1-m)>P(Y=1-y,M=m\mid X=m)$ holds, then there exists a non-empty subpopulation
$\Omega_{s}\subseteq\Omega$ such that for every $\omega\in\Omega_{s},$ $Y_{1}(\omega)=ym+(1-m)(1-y),$
$Y_{0}(\omega)=(1-y)m+(1-m)y,$ and $M_{1}(\omega)=M_{0}(\omega)=m.$
\end{theorem}

\begin{corollary}
Suppose $X$ is randomized at baseline. In addition suppose that $X$ has a positive monotonic effect on $M.$ If the inequality $P(Y=1,M=0\mid X=1)>P(Y=1,M=0\mid X=0)$ holds, then there exists a non-empty subpopulation $\Omega_{s}\subseteq \Omega$ such that for every $\omega\in\Omega_{s},$ the total effect, total direct effect, pure direct effect, principal stratum direct effect with $M_{1}(\omega)=M_{0}(\omega)=0$ are all equal to $1$, and the total indirect and pure indirect effects for $\omega$ are both zero.
\end{corollary}

\begin{proposition}
Suppose $X$ is randomized at baseline. In addition, suppose
that $X$ has a positive monotonic effect on $M.$ For some $m\in\{0,1\}$ and $y\in\{0,1\},$ the expression $P(Y=1-y,M=m\mid X=1-m)-P(Y=1-y,M=m\mid X=m)$ is equal to
\begin{align*}
 & P_c(ym+\{1-y\}\{1-m\},\{1-y\}m+\{1-m\}y,m,m) - P_c(m\{1-y\}+\{1-m\}y,ym+\{1-y\}\{1-m\},m,m)\\   &-P_c(m\{1-y\}+\{1-m\}y,ym+\{1-y\}\{1-m\},1,0) - P_c(1-y,1-y,1,0)
\end{align*}
\end{proposition}

\begin{corollary}
Suppose $X$ is randomized at baseline. In addition suppose that $X$ has a positive monotonic effect on $M.$ For $m\in\{0,1\}$ and $y\in\{0,1\}$, the expression $P(Y=1-y,M=m\mid X=1-m)-P(Y=1-y,M=m\mid X=m)$ is a lower bound on the counterfactual risk difference
\begin{align*}
&P_c(ym+\{1-y\}\{1-m\},\{1-y\}m+\{1-m\}y,m,m)-P_c(m\{1-y\}+\{1-m\}y,ym+\{1-y\}\{1-m\},m,m),&\nonumber
\end{align*}
and when $P(Y=1-y,M=m\mid X=1-m)>P(Y=1-y,M=m\mid X=m)$ the expression  $$\frac{P(Y=1-y,M=m\mid X=1-m)-P(Y=1-y,M=m\mid X=m)}{P(M=m\mid X=0)}\nonumber$$ is a lower bound on 
\begin{align}
\frac{
\left(\begin{array}{c}
P_c(ym+\{1-y\}\{1-m\},\{1-y\}m+\{1-m\}y,m,m)\\
-P_c(m\{1-y\}+\{1-m\}y,ym+\{1-y\}\{1-m\},m,m)
\end{array}\right)
}{P(M_1=M_0=m)},\nonumber
\end{align}
which relates to the standardized principal stratum direct effect. This standardized population principal stratum direct effect is equivalent to
\begin{align}
P(Y_1=ym+\{1-y\}\{1-m\},Y_0=\{1-y\}m+\{1-m\}y\mid M_1=m,M_0=m)\nonumber \\
-P(Y_1=m\{1-y\}+y\{1-m\},Y_0=ym+\{1-y\}\{1-m\}\mid M_1=m,M_0=m).\nonumber
\end{align}
\end{corollary}

{~~~~Proposition 7 enables similar interpretations of distribution of
the counterfactual contrasts that were provided for Proposition 2.
If~\(P(Y=1,M=0\mid X=1)>P(Y=1,M=0\mid X=0)\)}, then the proportion of
individuals~\(\omega\) randomized at baseline for whom treatment
has a beneficial effect on~\(Y\left(\omega\right)\) and~\(M_1\left(\omega\right)\)
and~\(M_0\left(\omega\right)\) are still zero is greater than the proportion of
individuals for whom treatment has a harmful effect
on~\(Y\) and~~\(M_1(\omega)\) and~\(M_0(\omega)\) are
still zero. The `instrumental inequalities' assuming~\(X\)
having a positive monotonic effect on~\(M\)
are~for~\(m\in\{0,1\}\) and
~\(y\in\{0,1\}\),~\(P(Y=1-y,M=m\mid X=1-m)\le P(Y=1-y,M=m\mid X=m)\) ~\citep{balke1997bounds}. These
`instrumental inequalities' serve as null hypotheses of the
corresponding principal stratum direct effect under investigation. At
most two of the four inequalities associated with Theorem 2 can hold.~

\begin{example}
Consider a population denoted $\Omega$ for which half of the individuals
$\omega_{1}$ follow counterfactual response type $Y_{1,M_{1}}(\omega_{1})=1,$
$Y_{0,M_{0}}(\omega_{1})=0,$ $M_{1}(\omega_{1})=0$ and $M_{0}(\omega_{1})=0$
and the other half of individuals $\omega_{2}$ follow counterfactual
response type $Y_{1,M_{1}}(\omega_{2})=0,$ $Y_{0,M_{0}}(\omega_{2})=1,$
$M_{1}(\omega_{2})=1$ and $M_{0}(\omega_{2})=1.$ We have $M_{1}(\omega)\geq M_{0}(\omega)$
for all $\omega\in\Omega$. For such a population, the average pure direct
effect, $E(Y_{1,M_{0}}-Y_{0,M_{0}}),$ and the total direct effect,
$E(Y_{1,M_{1}}-Y_{0,M_{1}}),$ are both zero. However, the two null
hypotheses associated with Theorem 2 are both falsified, indicating
that the two principal stratum direct effects associated with individual response types $\omega_1$ and $\omega_{2}$ are both present in our population. As noted in \citet{richardson2010analysis}, at most one of the the four inequalities associated with Theorem 1 can hold, which means that without the monotonicity assumption, one can at most detect one individual level principal stratum direct effect. However, two of the inequalities associated with Theorem 2 hold. We see that $P(Y=1,M=1\mid X=0)>P(Y=1,M=1\mid X=1)$ and $P(Y=1,M=0\mid X=1)>P(Y=1,M=0\mid X=0)$ in this example. This means that both individual level principal stratum direct effects can be detected using our results from Theorem 5.
\end{example}

~~~~Earlier, for the setting where one does not
assume~\(X\) has a positive monotonic effect
on~\(Y\), we saw Theorem 1 would quantify individual level
principal stratum and natural direct effects even when
the~\citet{sjolander2009bounds} bounds for the natural direct effect included
zero. The key difference between our results and~\citet{sjolander2009bounds}
results for the monotonicity setting is that we only assume
that~\(X\) has a monotonic effect on~\(M\),
wheras~\citet{sjolander2009bounds} makes additional monotonicity assumptions. We
are able to detect two different individual level principal stratum
direct effects even under the scenario where the total effect or natural
direct effect is zero. Detection of two different principal stratum
direct effects is not previously described in literature
~\citep{sjolander2009bounds,cai2008bounds,balke1997bounds,richardson2010analysis,Swanson_2018}.

\par\null

\subsection{Sensitivity analysis for monotonicity
assumptions}

{\label{473785}}

~~~~As demonstrated in Theorem 1, Corollary 2, and Proposition 3, we
provide methods to detect individual level direct effects and
characterize the magnitude of such individual level direct effects with
only assumptions that are guaranteed by design. To the best of our
knowledge, most methods assume that the nonparametric structural
equation model, influence diagram, or directed acyclic graph that is
used to point identify the direct effect is correctly specified and
available to the analyst~\citep{pearl2009causality}.~ Therefore, most if not all
approaches require assumptions that are not only unverifiable, but they
also offer no experimental design where these assumptions are guaranteed
by design.

~~~~Monotonicity assumptions are not guaranteed by design of experiment.
However, without these monotonicity assumptions, the tests to detect
individual level direct effects might be too stringent to detect even
large individual level direct effects in the population. In such
settings, where the investigator is unsure if the monotonicity
assumptions hold for all individuals in the population, informative
sensitivity analysis enable statisticians to interpret their results
even when these monotonicity assumptions are violated, and also
precisely quantify how violations of the monotonicity assumptions impact
the magnitude of the proposed individual level direct effect. The next
Theorem provides all the required ingredients necessary to conduct
informative and useful sensitivity analyses.

\begin{theorem}
Suppose $X$ is randomized at baseline. For $m\in\{0,1\}$ and $y\in\{0,1\},$ the expression $P(Y=1-y,M=m\mid X=1-m)-P(Y=1-y,M=m\mid X=m)-r$, where 
\begin{align*}
r=&P_c(ym+\{1-m\}\{1-y\},\{1-y\}m+\{1-m\}y,0,1)+P_c(1-y,1-y,0,1)
\nonumber\\
&-P_c(y\{1-m\}+\{1-y\}m,\{1-y\}\{1-m\}+ym,1,0)-P_c(1-y,1-y,1,0)&\nonumber
\end{align*}
is equal to the counterfactual risk difference
$$P_c(ym+\{1-m\}\{1-y\},\{1-y\}m+\{1-m\}y,m,m)-P_c(\{1-y\}m+y\{1-m\},\{1-y\}\{1-m\}+ym,m,m),\nonumber $$ 
and the expression
$$\frac{P(Y=1-y,M=m\mid X=1-m)-P(Y=1-y,M=m\mid X=m)-r}{P(M=m\mid X=1)-q}\nonumber,$$
where $q=P(M_1=m,M_0\neq m),$ is equal to
\begin{align*}
\frac{\left(\begin{array}{c}
P_c(ym+\{1-y\}\{1-m\},\{1-y\}m+\{1-m\}y,m,m)\\
-P_c(m\{1-y\}+\{1-m\}y,ym+\{1-y\}\{1-m\},m,m)
\end{array}\right)}{P(M_1=M_0=m)},
\end{align*}
which relates to the standardized principal stratum direct effect. This standardized population principal stratum direct effect is equivalent to
\begin{align*}
P(Y_1=ym+\{1-y\}\{1-m\},Y_0=\{1-y\}m+\{1-m\}y\mid M_1=m,M_0=m)\\
-P(Y_1=m\{1-y\}+y\{1-m\},Y_0=ym+\{1-y\}\{1-m\}\mid M_1=m,M_0=m).\nonumber
\end{align*}
\end{theorem}

~~~~The quantities~\(r\)
and~\(q\) are not point identified in randomized studies.
Consequently, the statistician or scientist will need to consider a range of suitable
values. A large~\(r\) will attenuate the the proposed
proportion of individual level direct effect as well as the standardized
principal stratum direct effect, while in contrast a
large~\(q\) will amplify the proposed standardized principal
stratum direct effect. As~\(r\) and~\(q\) are
functions only of potential outcomes, this enables scientists to
directly consider values of sensitivity parameters that are
intrinsically related to the physical experiment.
Here,~\(r\) should never exceed~\(P\left(Y=1-y,M=m\mid X=1-m\right)-P\left(Y=1-y,M=m\mid X=m\right)\)
and~\(q\) should never exceed~\(P\left(M=m\mid X=1\right).\)

~~~~The following example from a hypothetical randomized experiment
illustrates how sensitivity analysis parameter~~\(r\) can
be used to quantify how of~ individual level violations of the
monotonicity assumption can impact the studies findings. Suppose a
statistician estimates~\(P\left(Y=0,M=1\mid X=0\right)-P\left(Y=0,M=1\mid X=1\right)\) as 0.3. Then, for her study,
the statistician would need at least~\textbf{\(30\%\)} of the
individuals randomized in our study to follow counterfactual response
types~\(Y_1\left(\omega\right)=1,Y_{0\ }\left(\omega\right)=0,M_1\left(\omega\right)=0,M_0\left(\omega\right)=1\) or~\(Y_1\left(\omega\right)=0,Y_{0\ }\left(\omega\right)=0,M_1\left(\omega\right)=0,M_0\left(\omega\right)=1\) for there to be no
individual that follows the counterfactual response
type~\(Y_1\left(\omega\right)=1,Y_{0\ }\left(\omega\right)=0,M_1\left(\omega\right)=1,M_0\left(\omega\right)=1\) in her population. Similarly, if the scientist
reasons that at most~\(10\%\) of her study population follows
counterfactual response types~ ~\(Y_1\left(\omega\right)=1,Y_{0\ }\left(\omega\right)=0,M_1\left(\omega\right)=0,M_0\left(\omega\right)=1\)
or~\(Y_1\left(\omega\right)=0,Y_{0\ }\left(\omega\right)=0,M_1\left(\omega\right)=0,M_0\left(\omega\right)=1\), then she can conclude that at
least~\(20\%\) of her population follows counterfactual
response type~~\(Y_1\left(\omega\right)=1,Y_{0\ }\left(\omega\right)=0,M_1\left(\omega\right)=1,M_0\left(\omega\right)=1\). Prior studies or subject matter
expertise can help statisticians understand how individual level
violations of monotonicity assumptions can impact their findings.~

~~~~If the statistician is completely averse to sensitivity analysis,
then the results associated with Theorem 1, Corollary 2, Proposition 3,
and Corollary 4 are all still available to her. However, such
sensitivity analysis enables the statistician to detect and quantify
individual level directs when methods that assume randomization alone
are too stringent to detect and quantify individual level direct
effects. In this example, if the statistician takes~\(r=P\left(Y=1,M=0\mid X=1\right)+P\left(Y=0,M=0\mid X=1\right)\)
and~\(q=P\left(M=0\mid X=0\right),\) she will essentially be conducting the same
tests associated with Theorem 1 and Corollary 4 as expressions in
Theorem 9 with the proposed~\(q\) and~\(r\)
will take us back to exactly the same expressions as those provided in
Corollary 4.

~~~~Thus, the sensitivity analysis provided here can be viewed as a
compromise between the stringent setting where only baseline
randomization is assumed and the more flexible setting where in addition
to baseline randomization the statistician assumes positive
monotonicity. Some fixed~\(r\)~in Theorem 9 will provide
the true value of~\(P_c(y,1-y,m,m)-P_c(1-y,y,m,m)\) and for a
fixed~\(r\) and~\(q\) in Theorem 9 will
provide the standardized principal stratum direct effect, but
unfortunately there is no mechanism through which the statistician can
point identify this~\(r\)
and~\(q\).~\(\)Consequently, a range of
suitable values for~\(r\) and~\(q\) will
need to be considered, and scientists can deliberate whether the
proposed values are sensible.

\par\null

\section{Pleiotropic Effects}

{\label{132386}}

{~~~~In addition to variables and assumptions provided in Section 2,
let~\(Z\)~}denote an additional binary outcome. As before,
assume treatment~\(X\) is randomized, so
that~\(X\amalg (Z_{x},Y_{x}).\)~We also require the consistency assumption
for~\(Y_x\) and~\(Z_x\). This section will
define pleiotropic effects, which are of importance in the context of
making inferences abut the effect of a single gene or treatment on two
distinct outcomes or phenotypes~\citep{solovieff2013pleiotropy}. The results provided
herein enable scientists to discover such effects. Pleiotropic effects
are of interest in the detection of genetic variants that may operate
through multiple pathways or mechanisms. 

~~~~Pleiotropic effects are also of
interest in the context of Mendelian randomization studies that use a
genetic variant as an instrumental variable because pleiotropic effects
will often indicate that the exclusion restriction required for the
instrumental variable analysis is violated.~We give a formal
counterfactual definition for the intuitive notion of pleiotropic effects
as follows.

This section does not use notation~\(Y_{x,m}\)
or~\(Y_{x,M_{x^{\prime}}}\) from the previous three sections.

\begin{definition}[Pleiotropic Effects]
A treatment or exposure $X$ is defined to have an individual level
pleiotropic effect on outcomes $Y$ and $Z$ if there exists an individual
$\omega\in\Omega$ for whom one of the following four counterfactual
response types holds:
\begin{align*}
Y_{1}(\omega)=1,Y_{0}(\omega)=0,Z_{1}(\omega)=1,Z_{0}(\omega)=0;\nonumber\\
Y_{1}(\omega)=1,Y_{0}(\omega)=0,Z_{1}(\omega)=0,Z_{0}(\omega)=1;\nonumber\\
Y_{1}(\omega)=0,Y_{0}(\omega)=1,Z_{1}(\omega)=1,Z_{0}(\omega)=0;\nonumber\\
Y_{1}(\omega)=0,Y_{0}(\omega)=1,Z_{1}(\omega)=0,Z_{0}(\omega)=1.\nonumber
\end{align*}
\end{definition}

~~~~For the first condition, the treatment~\(X\)~is
causative for both~\(Y\)and~\(Z.\)~For the
second condition, the treatment~\(X\) is causative
for~\(Y\) and preventative for~\(Z.\) For the
third condition, the treatment is preventative for~\(Y\)
and causative for~\(Z.\) Finally, for the last condition,
\(X\)~is preventative for both~\(Z\)~and
\(Y.\)

\begin{theorem}
Suppose $X$ is randomized at baseline. If $P(Y=1,Z=1\mid X=1)+P(Y=0,Z=0\mid X=0)>1$, then there exists a non-empty subpopulation $\Omega_{s}\subseteq\Omega$ such that for all $\omega\in\Omega_{s},$ $Y_{1}(\omega)=1,$ $Y_{0}(\omega)=0$ and $Z_{1}(\omega)=1,$ $Z_{0}(\omega)=0.$ \\

The expression $P(Y=1,Z=1\mid X=1)+P(Y=0,Z=0\mid X=0)-1$ is equal to $$P(Y_1=1,Y_0=0,Z_1=1,Z_0=0)-\sum_{(y_1,z_1)\in\{0,1\}^2 \setminus \{1,1\}}P(Y_1=y_1,Z_1=z_1,(Y_0,Z_0)\neq (0,0)),\nonumber$$ and is therefore a lower bound for $P(Y_1=1,Y_0=0,Z_1=1,Z_0=0).$
\end{theorem}

~~~~The other three counterfactual response types can be likewise
detected using analogous results collected in the online supplement.
Theorem 10 enables the detection of individual level pleiotropic
effects, but does not distinguish between two different forms of
individual level pleiotropy. The first situation is
where~\(X\) causes~\(Y\) which in turn
causes~\(Z\) (or~\(X\)
causes~\(Z\) which in turn causes~\(Y\)). The
second situation is where ~\(X\)causes~\(Y\)
through a pathway not through~\(Z\)~ and~\(X\)
causes~\(Z\) through a pathway not
through~\(Y\). The first situation is called `mediated
pleiotropy' and the second situation is known as `biologic
pleiotropy'~\citep{solovieff2013pleiotropy}. Theorem 10 allows an investigator to
detect if at least one of these two types of pleiotropy is present. ~

\par\null

\section{Data Analysis}

{\label{198257}}

~~~~\citet{yerushalmy1964mother}
noted that maternal smoking~\(X\in\left\{0,1\right\}\) is associated with lower
infant mortality~\(Y\in\left\{0,1\right\}\) among low birth weight
infants~\(M\in\left\{0,1\right\}\). Other authors after Yerushalmy~\citep{hernandez2006birth,vanderweele2012conditioning} have replicated this finding, and referred to this finding as the birth weight paradox, as the direction of the association is
opposite of what is expected. Here,~\(X=1\) denotes that the
mother smoked cigarettes, and~\(X=0\) denotes that the
mother did not smoke cigarettes;~\(M=1\) denotes that the
mother had a low birth weight baby (less than 2500 grams),
and~\(M=0\) denotes that the mother did not have a low birth
weight baby;~\(Y=1\) denotes that the the baby died one year
after birth, and~\(Y=0\) denotes that the baby lived at
least one year after birth. Specifically,~\citet{yerushalmy1964mother}~reported
finding that, among low birth weight infants of less than 2500 grams,
the observed risk of infant mortality was greater for non-smoking
mothers than non-smoking mothers. The `paradox' stems from the
realization that smoking is known to be detrimental to human health, but
in the analysis of~\citet{yerushalmy1964mother} it seemingly appeared as though
for babies that are born with low birth weight maternal smoking was
protective for infant mortality. The two sets of 2x2 contingency tables
are provided to ensure that our results are reproducible.~~
\par\null\par\null\selectlanguage{english}
\begin{table}[h!]
{\begin{center}
\caption{{Yerushalmy's Infant Mortality Data}
{\label{109891}}}%
\normalsize\begin{tabulary}{1.0\textwidth}{CCCCCCCC}
 & X=1 &  &  &  & X=0 &  &  \\
 & Y=1 & Y=0 & Total &  & Y=1 & Y=0 & Total \\
M=1 & 27 & 210 & 237 & M=1 & 43 & 154 & 197 \\
M=0 & 15 & 3474 & 3489 & M=0 & 24 & 5846 & 5870 \\
Total & 42 & 3684 & 3726 & Total & 67 & 6000 & 6067 \\
\end{tabulary}
\end{center}
}
\end{table}

~~~~Table 3 in~\citet{yerushalmy2014relationship} provides us with the raw numbers for
the first row of the multi-way contingency table provided in Table 1
below. The second row of Table 1 comes from using Table 2
of~\citet{yerushalmy2014relationship} to calculate the raw numbers of neonatal deaths
for smokers and nonsmokers from the mortality rates. The column totals
are provided by~\citep{yerushalmy2014relationship} as well as the cell entries for the
very first row. With the mortality rates for smokers and nonsmokers, we
are able to recreate the raw numbers of infant mortality for normal
birth weight babies.

~~~~~In Yerushalmy's data of white mothers, the risk
difference of infant mortality comparing low birth weight infants whose
mothers smoke versus those who do not is~\(P\left(Y=1\mid X=1,M=1\right)-P\left(Y=1\mid X=0,M=1\right)=27/237-43/197=0.114-0.218=-0.104\) with a
two-sided confidence interval of~~\(\left(-0.18,-0.03\right)\) \citep{yerushalmy2014relationship}.~ This result might
erroneously, foolishly and catastrophically interpreted as evidence that
maternal smoking reduces the risk of infant mortality for low birth
weight babies~\citep{yerushalmy1964mother}. However, if one posits the existence
of a nonparametric structural equation model that generates the observed
data~\citep{pearl2009causality}, such analysis ignores the possibility of
unmeasured factors that are common causes of both low birth
weight~\(M\) and infant mortality~\(Y\).
Previous research speculated that malnutrition or birth defects could
act as common causes of low birth weight and infant
mortality~~\citep{hernandez2006birth,vanderweele2012conditioning}.~

~~~~Even if the relationship between maternal smoking and birth weight
is not confounded, and also that between maternal smoking and infant
mortality, an unmeasured common cause of low birth weight and infant
mortality could introduce bias of the analysis. The results provided in
our paper do not require having data on, or having controlled for,
potential common causes of low birth weight and infant mortality.
Furthermore, our analysis is agnostic to whether the data is generated
by nonparametric structural equation models, influence diagrams or
directed acyclic graphs~\citep{pearl2009causality} with untestable and
unverifiable assumptions, which is required for prior research
concerning the birth weight paradox.

~~~~We will assume that the population of white mothers that smoke are
exchangeable with mothers that do not smoke, which is
denoted~\(\left(Y_x,M_x\right)\amalg X\). We will also assume that there is no baby
that would be born with normal birth weight had the mother smoked, but
this same baby would be born with low birth weight had the mother not
smoked, which is the positive monotonicity
assumption~\(M_1\left(\omega\right)\ge M_0\left(\omega\right)\)~ for all~\(\omega\in\Omega.\)~ Using
Yerushalmy's data, we find that the risk of giving birth to a low birth
weight baby is nearly double for smokers as compared to non-smoking
mothers. Vast literature on reproductive, perinatal and pediatric
epidemiology provides strong evidence that smoking is likely to cause
low birth weight babies~\citep{Hebel1988,Hebel1985,Sexton1984}.

~~~~We can now apply Theorem 5, Proposition 7 and Corollary 8. To detect
whether there exists at least one baby that would be born with low birth
weight regardless of the mother's smoking status and maternal smoking
prevents infant mortality for the same baby, we would attempt to falsify
the null hypothesis~\(P(Y=1,M=1\mid X=0)\le P(Y=1,M=1\mid X=1)\). If we falsify this null
hypothesis, then from Corollary 8, we can conclude that the proportion
of babies that are born with low birth weight regardless of maternal
smoking status and smoking protects the baby from infant mortality is
greater than the proportion of babies that are born with low birth
weight regardless of maternal smoking status and smoking causes the
baby's death within the first year of birth. From Proposition 7, the
risk difference~~\(P(Y=1,M=1\mid X=0)-P(Y=1,M=1\mid X=1)\) is equal to~\(P_c\left(0,1,1,1\right)-\left\{P_{c\ }\left(1,0,1,1\right)+P_c\left(1,0,1,0\right)+P_c\left(1,1,1,0\right)\right\}\).
Therefore, the same risk difference is a lower bound
on~\(P_{c\ }\left(0,1,1,1\right)-P_c\left(1,0,1,1\right)\) from Corollary 8.

~~~~From Corollary 8, the lower bound on the proportion of babies that
would be born with low birth weight regardless of maternal smoking and
these babies would be prevented from dying in its first year from
maternal smoking is given as~\(P(Y=1,M=1\mid X=0)-P(Y=1,M=1\mid X=1).\)~ Using Yerushalmy's
data this risk difference is ~\(43/6067-27/3726=0.0071-0.0072=-0.0001\)~ with a
one-sided~\(\)95 percent confidence interval
of~\((-0.003,1)\). We therefore do not find evidence that smoking
prevents infant mortality for babies that would be born low birth weight
regardless of maternal smoking status. Had we found evidence for such
babies in our population using Corollary 8, we could have applied
Theorem 9 to assess the consequences of the monotonicity assumptions on
exaggerating the proportion of babies that regardless of maternal
smoking status are born with low birth weight and maternal smoking
prevents infant mortality.

~~~~To test the more reasonable hypothesis, given the preponderance of
scientific evidence of the harmful effect of smoking on human health,
that smoking actually causes infant mortality for babies that would be
born low birth weight regardless of maternal smoking, we again can apply
Theorem 5, Proposition 7 and Corollary 8.~ From Theorem 5, to detect
whether there exists at least one baby that would be born with low birth
weight regardless of maternal smoking and maternal smoking causes infant
mortality for the same baby, we would attempt to falsify the null
hypothesis~\(P(Y=0,M=1\mid X=0)\le P(Y=0,M=1\mid X=1)\). If we falsify this null hypothesis, then
from Corollary 8, we can conclude that the proportion of babies that are
born with low birth weight regardless of maternal smoking status and
smoking causes the baby's infant mortality is greater than the
proportion of babies that are born with low birth weight regardless of
maternal smoking status and smoking prevents the baby's infant mortality
within the first year of birth.~

~~~~From Proposition 7, the risk difference~~\(P(Y=0,M=1\mid X=0)-P(Y=0,M=1\mid X=1)\) is
equal to~\(P_c\left(1,0,1,1\right)-\left\{P_{c\ }\left(0,1,1,1\right)+P_c\left(0,1,1,0\right)+P_c\left(0,0,1,0\right)\right\}\). From Corollary 8, the same risk
difference~ is also a lower bound on~\(P_{c\ }\left(1,0,1,1\right)-P_c\left(0,1,1,1\right)\). ~From Theorem
5 and Corollary 8, the lower bound on the proportion for a baby that
would be born low birth weight regardless of maternal smoking and this
baby would die in its first year from maternal smoking is given
as~\(P(Y=0,M=1\mid X=0)-P(Y=0,M=1\mid X=1).\)~ Using Yerushalmy's data this risk difference is
~\(154/6067-210/3726=0.025-0.056=-0.031\)~ with a one-sided~\(\)95 percent
confidence interval of~\((-0.038,1)\). Therefore, using Corollary
8, we do not find evidence that smoking causes infant mortality for
babies that would be born low birth weight regardless of maternal
smoking status. Had we found evidence for such babies in our population
using Corollary 8, we could have applied Theorem 9 to assess the
consequences of the monotonicity assumptions on exaggerating the
proportion of babies that regardless of maternal smoking status are born
with low birth weight and maternal smoking causes infant mortality.

~~~~To understand why we are unable to find statistical evidence for
this second more reasonable hypothesis, we can conduct a sensitivity
analysis similar to what is presented in Theorem 9. As stated earlier,
under the reasonable monotonicity assumption that smoking never prevents
low birth weight, the risk difference~~\(P(Y=0,M=1\mid X=0)-P(Y=0,M=1\mid X=1)\) is equal
to~\(P_c\left(1,0,1,1\right)-\left\{P_{c\ }\left(0,1,1,1\right)+P_c\left(0,1,1,0\right)+P_c\left(0,0,1,0\right)\right\}\).~ ~ ~This counterfactual risk difference includes
the term~\(P_c\left(0,0,1,0\right)\), which is the proportion of babies that as
a result of smoking would not die in their first year after birth
regardless of maternal smoking, but would be born of low birth weight
only if the mother smoked. From yet another application of Corollary 8
(interchanging the roles of \(M\) and
\(Y\)), under the reasonable assumption that smoking never
prevents infant mortality, a lower bound on this
proportion~\(P_c\left(0,0,1,0\right)\) is given by~\(P\left(M=1,Y=0\mid X=1\right)-P\left(M=1,Y=0\mid X=0\right)=210/3726-154/6067=0.031\) which has
a one-sided 95 percent confidence interval of~~\(\left(0.024,1\right)\),
suggesting that~~\(P_c\left(0,0,1,0\right)\) is at least~\(0.24\).

~~~~Under the monotonicity assumption that smoking never prevents a low
birth weight baby, the risk difference~\(P(Y=0,M=1\mid X=0)+P_c\left(0,0,1,0\right)-P(Y=0,M=1\mid X=1)\)~ is equal to the counterfactual expression~\(P_c\left(1,0,1,1\right)-\left\{P_{c\ }\left(0,1,1,1\right)+P_c\left(0,1,1,0\right)\right\}\). While we cannot point identify the risk
difference~\(P(Y=0,M=1\mid X=0)+P_c\left(0,0,1,0\right)-P(Y=0,M=1\mid X=1),\) it appears as though it is non-negative
in Yerushalmy's data. Therefore, Yerushalmy's data contains no evidence
that smoking prevents infant mortality for babies that regardless of
maternal smoking would be born with low birth weight from the very first
hypothesis test that we conducted. In fact, under realistic monotonicity
assumptions and reasonable sensitivity analyses, the evidence is
suggestive that maternal smoking actually causes infant mortality for
babies that regardless of maternal smoking would be born with low birth
weight. Hypothesis tests to detect whether regardless of maternal
smoking a baby would be born with normal birth weight and smoking causes
infant mortality could likewise be conducted using our methods.

\par\null

\section{Generalizations and
Extensions}

{\label{465172}}

\subsection{Multivariate Direct Effects with Continuous and Discrete
Outcomes and
Mediators}
{\label{334479}}
Our previous results are derived for a binary
outcome~\(Y\) and~ binary mediator~\(M\). We
now provide results for tuple
valued~\(\)\(Y\) and~\(M,\)
each of which has support on~\(\mathbf{R}^{n_y\ }\) and~\(\mathbf{R}^{n_m},\)
where~\(n_y\) and~\(n_m\) denote the number of
elements in the corresponding tuple. Here,~\(Y\)
and~\(M\) could have continuous and/or discrete
components. The special case of only continuous components is provided
in the online supplement. All of our previous results can be interpreted
as specializations of the results presented in this section.
Treatment~\(X\) is still binary, but categorical or
ordinal treatments can be suitably dichotomized to apply the results in
this section at any level of the dichotomization. In this section, we
allow~\(Y_1\left(\omega\right)\),~~\(Y_0\left(\omega\right),\)~\(M_1\left(\omega\right)\)
and~\(M_0\left(\omega\right)\) to be stochastic counterfactuals that have
distribution function~\(F_{Y_1,Y_0,M_1,M_0}\left(y_a,y_b,m_a,m_b\right)=P\left(Y_1\in y_a,Y_0\in y_b,M_1\in m_a,M_0\in m_b\right)\),~ where~\(y_a\)
and~\(y_b\)~are subsets of~\(\mathbf{R}^{n_y}\)
and~\(m_a\) and~\(m_b\) are subsets
of~\(\mathbf{R}^{n_m}\). We still require the randomization assumption
that~\(\left(Y_x,M_{x\ }\right)\) is independent of~\(X\) and the
consistency assumption~ that
when~\(X\left(\omega\ \right)=x\),~\(Y_x\left(\omega\right)=Y\left(\omega\right)\) and~\(M_x\left(\omega\right)=M\left(\omega\right)\).
\begin{proposition}
Suppose $X$ is randomized at baseline. For some $y_{a}\subset\mathbf{R}^{n_{y}}$ and $m_{a}\mathbf{\subset R}^{n_{m}}$, the contrast  $$P(Y\in y_{a},M\in m_{a}\mid X=1)+P(Y\not\in y_{a},M\in m_{a}\mid X=0)-1\nonumber$$
 is equal to
\begin{align*}
 & P(Y_{1}\in y_{a},Y_{0}\not\in y_{a},M_{1}\in m_{a},M_{0}\in m_{a})-P\left(Y_{1}\not\in y_{a},Y_{0}\in y_{a},M_{1}\in m_{a},M_{0}\in m_{a}\right)\nonumber\\
 & -\left\{ P\left(Y_{1}\in y_{a},Y_{0}\in y_{a},M_{1}\not\in m_{a},M_{0}\in m_{a}\right)+P\left(Y_{1}\not\in y_{a},Y_{0}\in y_{a},M_{1}\not\in m_{a},M_{0}\in m_{a}\right)\right.\nonumber\\
 & +P\left(Y_{1}\in y_{a},Y_{0}\not\in y_{a},M_{1}\in m_{a},M_{0}\not\in m_{a}\right)+P\left(Y_{1}\not\in y_{a},Y_{0}\not\in y_{a},M_{1}\in m_{a},M_{0}\not\in m_{a}\right)\nonumber\\
 & +P\left(Y_{1}\in y_{a},Y_{0}\not\in y_{a},M_{1}\not\in m_{a},M_{0}\not\in m_{a}\right)+P\left(Y_{1}\not\in y_{a},Y_{0}\not\in y_{a},M_{1}\not\in m_{a},M_{0}\not\in m_{a}\right)\nonumber\\
 & \left.+P\left(Y_{1}\in y_{a},Y_{0}\in y_{a},M_{1}\not\in m_{a},M_{0}\not\in m_{a}\right)+P\left(Y_{1}\not\in y_{a},Y_{0}\in y_{a},M_{1}\not\in m_{a},M_{0}\not\in m_{a}\right)\right\},\nonumber
\end{align*}
and therefore is a sharp lower bound of $P(Y_{1}\in y_{a},Y_{0}\not\in y_{a},M_{1}\in m_{a},M_{0}\in m_{a})-P(Y_{1}\not\in y_{a},Y_{0}\in y_{a},M_{1}\in m_{a},M_{0}\in m_{a}).$

When  $P(Y\in y_{a},M\in m_{a}\mid X=1)+P(Y\not\in y_{a},M\in m_{a}\mid X=0)>1,\nonumber$ the expression
$$\frac{P(Y\in y_{a},M\in m_{a}\mid X=1)+P(Y\not\in y_{a},M\in m_{a}\mid X=0)-1}{P(M\in m_{a}\mid X=1)-P(M\not\in m_{a}\mid X=0)}\nonumber$$ is a sharp lower bound on the standardized principal stratum direct effect
$$\frac{P(Y_{1}\in y_{a},Y_{0}\not\in y_{a},M_{1}\in m_{a},M_{0}\in m_{a})-P\left(Y_{1}\not\in y_{a},Y_{0}\in y_{a},M_{1}\in m_{a},M_{0}\in m_{a}\right)}{P(M_{1}\in m_{a},M_{0}=m_{a})},\nonumber$$
which is equivalent to 
$$P(Y_{1}\in y_{a},Y_{0}\not\in y_{a}\mid M_{1}\in m_{a},M_{0}\in m_{a})-P\left(Y_{1}\not\in y_{a},Y_{0}\in y_{a}\mid M_{1}\in m_{a},M_{0}\in m_{a}\right).\nonumber$$
\end{proposition}

This demonstrates that our results apply are also applicable to
continuous random variables. This work demonstrates that quantification
of direct effects can be conducted with only design based assumptions.
Such design based assumptions make the study of direct effects possible
with randomized trials. Consider now the analogous results under
monotonicity assumptions.

\begin{proposition}
Suppose $X$ is randomized at baseline. In addition, suppose there exists no individual $\omega\in\Omega$ of response type $M_1(\omega)\not\in m_{a}$ and $M_0(\omega)\in m_a$.  For some $y_{a}\subset\mathbf{R}^{n_{y}}$ and $m_{a}\mathbf{\subset R}^{n_{m}}$, the contrast  $$P(Y\not\in y_{a},M\in m_{a}\mid X=0)-P(Y\not\in y_{a},M\in m_{a}\mid X=1)\nonumber$$ is a sharp lower bound on
$$P(Y_{1}\in y_{a},Y_{0}\not\in y_{a},M_{1}\in m_{a},M_{0}\in m_{a})
-P(Y\not\in y_{a},Y_{0}\in y_{a},M_{1}\in m_{a},M_{0}\in m_{a}).\nonumber$$
When there exists no individual $\omega\in\Omega$ of response type $M_1(\omega)\not\in m_{a}$ and $M_0(\omega)\in m_a$ and 
$$P(Y\not\in y_{a},M\in m_{a}\mid X=0)-P(Y\not\in y_{a},M\in m_{a}\mid X=1)>0\nonumber$$
then the empirical condition 
$$\frac{P(Y\not\in y_{a},M\in m_{a}\mid X=0)-P(Y\not\in y_{a},M\in m_{a}\mid X=1)}{P(M\in m_{a}\mid X=0)}\nonumber$$
serves as a sharp lower bound on the standardized principal stratum direct effect
$$\frac{P(Y_{1}\in y_{a},Y_{0}\not\in y_{a},M_{1}\in m_{a},M_{0}\in m_{a})-P(Y\not\in y_{a},Y_{0}\in y_{a},M_{1}\in m_{a},M_{0}\in m_{a})}{P(M_1\in m_a,M_0\in m_a)},\nonumber$$ which is equivalent to 
$$P(Y_{1}\in y_{a},Y_{0}\not\in y_{a}\mid M_{1}\in m_{a},M_{0}\in m_{a})-P(Y\not\in y_{a},Y_{0}\in y_{a}\mid M_{1}\in m_{a},M_{0}\in m_{a}).\nonumber$$ 
\end{proposition}

For simple sensitivity analysis for the monotonicity assumptions, we
present the following results. These results enable scientists to
directly consider how violations of monotonicity assumptions impact the
findings of their results. Monotonicity assumptions enable scientists to
uncover unit level direct effects in cases where without such
assumptions the resulting tests are too stringent. Consequently,
informative and simple sensitivity analyses that enables scientists to
understand how such assumptions impact study results is important.

\begin{proposition}
Suppose $X$ is randomized at baseline. The expression $P(Y\not\in y_{a},M\in m_{a}\mid X=0)-P(Y\not\in y_{a},M\in m_{a}\mid X=1)-r$, where
\begin{align*}
r=& \left\{ P(Y_{1}\not\in y_{a},Y_{0}\not\in y_{a},M_{1}\not\in m_{a},M_{0}\in m_{a})+P(Y_{1}\in y_{a},Y_{a}\not\in y_{a},M_{1}\not\in m_{a},M_{0}\in m_{a})\nonumber\right.\\
 & -\left.P(Y_{1}\not\in y_{a},Y_{0}\not\in y_{a},M_{1}\in m_{a},M_{0}\not\in m_{a})-P(Y_{1}\not\in y_{a},Y_{0}\in y_{a},M_{1}\in m_{a},M_{0}\not\in m_{a})\right\},\nonumber\nonumber
\end{align*}
is equal to the unstandardized principal stratum direct effect, $P(Y_{1}\in y_{a},Y_{0}\not\in y_{a},M_{1}\in m_{a},M_{0}\in m_{a})-P(Y\not\in y_{a},Y_{0}\in y_{a},M_{1}\in m_{a},M_{0}\in m_{a}).$
Additionally, the expression
$$\frac{P(Y\not\in y_{a},M\in m_{a}\mid X=0)-P(Y\not\in y_{a},M\in m_{a}\mid X=1)-r}{P(M\in m_{a}\mid X=0)-q},\nonumber$$
where $q=P(M_{1}\not\in m_{a},M_{0}\in m_{a}),$ is equal to the standardized principal stratum direct effect $$\frac{P(Y_{1}\in y_{a},Y_{0}\not\in y_{a},M_{1}\in m_{a},M_{0}\in m_{a})-P(Y\not\in y_{a},Y_{0}\in y_{a},M_{1}\in m_{a},M_{0}\in m_{a})}{P(M_{1}\in m_{a},M_{0}\in m_{a})}, \nonumber$$
which is equivalent to 
$$P(Y_{1}\in y_{a},Y_{0}\not\in y_{a}\mid M_{1}\in m_{a},M_{0}\in m_{a})-P(Y\not\in y_{a},Y_{0}\in y_{a}\mid M_{1}\in m_{a},M_{0}\in m_{a})\nonumber.$$
\end{proposition}~~~Note our results in previous sections with binary~\(Y\)
and~\(M\) are specializations of the results provided here
with special sets for~\(y_a\) and~\(m_a\). The
full generality of our results provided here are appealing, but the
specificity of our previous results are easier to read and comprehend.
Under randomization, the results provided in this section demonstrate
that coarsening a continuous or ordinal outcome or mediator and still
provide useful results.

\subsection{Direct Effects Through Time With and Without Missing
Data~}

{\label{734122}}

~~~~~When there is no missing data, instead of
examining~\(Y\) and~\(M\) at one specific time
after randomization, the scientist can consider stochastic
processes~\(\left(Y_{\ }\left(\omega,t\right),M_{ }\left(\omega,t_1\right)\right)\) for any set of
times~\(t\)~and~\(t_1\le t\) after randomization, and
apply any or all of our results for each specific~\(t\)
and~\(t_1.\) All the statistician has to do is
replace~\(Y\left(\omega\right)\)
and~\(M\left(\omega\right)\)~with~\(Y\left(\omega,t\right)\) and~\(M\left(\omega,t_1\right)\)
respectively and~\(Y_x\left(\omega\right)\) and~\(M_x\left(\omega\right)\)
with~\(Y_x\left(\omega,t\right)\) and~\(M_x\left(\omega,t_1\right)\) in all our Theorems,
Propositions and Corollaries. We provide the full set of results in the
online supplement for completeness.~

~~~~ If the statistician has missing data or censoring in the dataset,
then this can be easily handled. The tuple
valued~\(M_{ }\left(\omega,t_1\right)\)~for each
specific~\(t\)~and~\(t_1\) can include a
variable~\(R\left(\omega,t_1\right),\) called response indicators, that takes
value 0 if either~\(Y\) or~\(M\)~is missing
for individual~\(\omega\) at
time~\(t\)~for~\(y\) or~\(t_1\)
for~\(M,\) and 1 otherwise for
each~individual~\(\omega\) at time pair~\(\left(t,t_1\right)\) ,
and then apply Proposition 11, 12 or 14 as required. For more refined
handling of the missing data as well as sensitivity analysis, we refer
the reader to the online supplement of this paper which provides methods
to handle missing data for the mediator and/or outcome with and without
monotonicity assumptions. 

~~~~The time index can be tuple valued as well, or $t=(t_1,\ldots,t_k).$ In this setting, $Y(t)$ would denote the trajectories $Y(t)=(Y(t_1),\ldots,Y(t_k))$ and $M(t)=(M(t_1),\ldots,M(t_k)).$ The counterfactuals are similarly $Y_x(t)=(Y_x(t_1),\ldots,Y_x(t_k)),$ and $M_x(t)=(M(t_1),\ldots,M_x(t_k)).$ Substituting these factuals and counterfactuals in our results above, scientists are able to examine individual level direct effects through time. If $M(t)$ and $M_x(t)$ are taken to be the empty set, then we return to the setting of total effects for $Y(t)$ through time.

\section{Conclusion}

{\label{476172}}

This work demonstrates that quantifying individual level natural direct
effects, principal stratum direct effects, and pleiotropic effects can
be implemented with no stronger identification assumptions than testing
for a non-zero total effect. The identifiability assumptions that are
needed to test the efficacy of a drug in a randomized clinical trial are
exactly the identifiability assumptions needed to test for the relevant
individual level direct effects. This embeds testing for these direct
effects firmly within the Neyman-Pearson paradigm. Previous literature
studied the `instrumental inequalities' and made contributions to
understanding some of their consequences, including that falsification
of the `instrumental inequality' provides evidence for a non-zero
population average controlled direct effect~\citep{cai2008bounds}. Our
paper is the first to establish~tests for an individual level principal
stratum direct effects and natural direct effects without unverifiable
and at times uninterpretable cross world or graphical assumptions. We
also quantify the magnitude of such direct effects. Under monotonicity,
our conditions can be used to detect individual level natural direct
effects even when the bounds for the population average natural direct
effects include zero. We also generate results that provide information
on the magnitude of the corresponding principal stratum direct effects.
We enable statisticians to conduct sensitivity analysis of monotonicity
assumptions to yield interpretable results for scientific investigation.
Our methodology can handle all types of outcomes and mediators: binary,
categorical, continuous, mixed, multivariate, and multivariate
measurements with missing data that are stochastic processes through
time. Additionally, we define a new causal effect, individual level
pleiotropy, in the counterfactual or potential outcome framework, and
derive the associated empirical conditions which could be used to detect
such an effect in a population.

\section*{Online Supplement}

{\label{875488}}

Tables that provide further aid to readers to derive our results are
provided in the online supplement. We also discuss the case of
continuous, ordinal, multivariate, time series mediators and outcomes in
the online supplement, as well as more detailed proofs for tuple valued
outcomes and mediators with stochastic counterfactuals.~

\par\null

\section*{\texorpdfstring{{Acknowledgements}}{Acknowledgements}}

{\label{362067}}

{We thank James M. Robins, Eric J. Tchetgen Tchetgen, Michael Hudgens,
and Linbo Wang for helpful discussions. The authors also acknowledge
useful comments from two anonymous reviewers and the associate editor.
This research was supported by the National Institutes of Health.}

\selectlanguage{english}

\bibliography{converted_to_latex.bib%
}

\begin{thebibliography}{25}
\providecommand{\natexlab}[1]{#1}
\providecommand{\url}[1]{\texttt{#1}}
\expandafter\ifx\csname urlstyle\endcsname\relax
  \providecommand{\doi}[1]{doi: #1}\else
  \providecommand{\doi}{doi: \begingroup \urlstyle{rm}\Url}\fi

\bibitem[Avin et~al.(2005)Avin, Shpitser, and Pearl]{avin2005identifiability}
Chen Avin, Ilya Shpitser, and Judea Pearl.
\newblock Identifiability of path-specific effects.
\newblock In \emph{Proceedings of the 19th international joint conference on
  Artificial intelligence}, pages 357--363. Morgan Kaufmann Publishers Inc.,
  2005.

\bibitem[Balke and Pearl(1997)]{balke1997bounds}
Alexander Balke and Judea Pearl.
\newblock Bounds on treatment effects from studies with imperfect compliance.
\newblock \emph{Journal of the American Statistical Association}, 92\penalty0
  (439):\penalty0 1171--1176, 1997.

\bibitem[Cai et~al.(2008)Cai, Kuroki, Pearl, and Tian]{cai2008bounds}
Zhihong Cai, Manabu Kuroki, Judea Pearl, and Jin Tian.
\newblock Bounds on direct effects in the presence of confounded intermediate
  variables.
\newblock \emph{Biometrics}, 64\penalty0 (3):\penalty0 695--701, 2008.

\bibitem[Frangakis and Rubin(2002)]{frangakis2002principal}
Constantine~E Frangakis and Donald~B Rubin.
\newblock Principal stratification in causal inference.
\newblock \emph{Biometrics}, 58\penalty0 (1):\penalty0 21--29, 2002.

\bibitem[Hebel et~al.(1985)Hebel, Nowicki, and Sexton]{Hebel1985}
J~Richard Hebel, Patricia Nowicki, and Mary Sexton.
\newblock The effect of antismoking intervention during pregnancy: an
  assessment of interactions with maternal characteristics.
\newblock \emph{Am J Epidemiol}, 122:\penalty0 135--48, Jul 1985.

\bibitem[Hebel et~al.(1988)Hebel, Fox, and Sexton]{Hebel1988}
J~Richard Hebel, Norma~Lynn Fox, and Mary Sexton.
\newblock Dose-response of birth weight to various measures of maternal smoking
  during pregnancy.
\newblock \emph{J Clin Epidemiol}, 41:\penalty0 483--9, 1988.

\bibitem[Hern{\'a}ndez-D{\'\i}az et~al.(2006)Hern{\'a}ndez-D{\'\i}az,
  Schisterman, and Hern{\'a}n]{hernandez2006birth}
Sonia Hern{\'a}ndez-D{\'\i}az, Enrique~F Schisterman, and Miguel~A Hern{\'a}n.
\newblock The birth weight ``paradox uncovered?''.
\newblock \emph{American journal of epidemiology}, 164\penalty0 (11):\penalty0
  1115--1120, 2006.

\bibitem[Hudgens et~al.(2003)Hudgens, Hoering, and Self]{Hudgens2003}
Michael~G Hudgens, Antje Hoering, and Steven~G Self.
\newblock On the analysis of viral load endpoints in hiv vaccine trials.
\newblock \emph{Stat Med}, 22:\penalty0 2281--98, Jul 2003.

\bibitem[Imai et~al.(2013)Imai, Tingley, and Yamamoto]{imai2013experimental}
Kosuke Imai, Dustin Tingley, and Teppei Yamamoto.
\newblock Experimental designs for identifying causal mechanisms.
\newblock \emph{Journal of the Royal Statistical Society: Series A (Statistics
  in Society)}, 176\penalty0 (1):\penalty0 5--51, 2013.

\bibitem[Imbens and Rubin(2015)]{imbens2015causal}
Guido~W Imbens and Donald~B Rubin.
\newblock \emph{Causal inference in statistics, social, and biomedical
  sciences}.
\newblock Cambridge University Press, 2015.

\bibitem[Pearl(2009)]{pearl2009causality}
Judea Pearl.
\newblock \emph{Causality}.
\newblock Cambridge University Press, 2009.

\bibitem[Ramsahai(2013)]{ramsahai2013probabilistic}
Roland~R Ramsahai.
\newblock Probabilistic causality and detecting collections of interdependence
  patterns.
\newblock \emph{Journal of the Royal Statistical Society: Series B (Statistical
  Methodology)}, 75\penalty0 (4):\penalty0 705--723, 2013.

\bibitem[Richardson and Robins(2010)]{richardson2010analysis}
Thomas~S Richardson and James~M Robins.
\newblock Analysis of the binary instrumental variable model.
\newblock \emph{Heuristics, Probability and Causality: A Tribute to Judea
  Pearl}, pages 415--444, 2010.

\bibitem[Robins and Greenland(1992)]{robins1992identifiability}
James~M Robins and Sander Greenland.
\newblock Identifiability and exchangeability for direct and indirect effects.
\newblock \emph{Epidemiology}, pages 143--155, 1992.

\bibitem[Sexton and Hebel(1984)]{Sexton1984}
Mary Sexton and J.~Richard Hebel.
\newblock A clinical trial of change in maternal smoking and its effect on
  birth weight.
\newblock \emph{JAMA}, 251:\penalty0 911--5, Feb 1984.

\bibitem[Sj{\"o}lander(2009)]{sjolander2009bounds}
Arvid Sj{\"o}lander.
\newblock Bounds on natural direct effects in the presence of confounded
  intermediate variables.
\newblock \emph{Statistics in Medicine}, 28\penalty0 (4):\penalty0 558--571,
  2009.

\bibitem[Solovieff et~al.(2013)Solovieff, Cotsapas, Lee, Purcell, and
  Smoller]{solovieff2013pleiotropy}
Nadia Solovieff, Chris Cotsapas, Phil~H Lee, Shaun~M Purcell, and Jordan~W
  Smoller.
\newblock Pleiotropy in complex traits: challenges and strategies.
\newblock \emph{Nature Reviews Genetics}, 14\penalty0 (7):\penalty0 483--495,
  2013.

\bibitem[Swanson et~al.(2018)Swanson, Hern{\'{a}}n, Miller, Robins, and
  Richardson]{Swanson_2018}
Sonja~A. Swanson, Miguel~A. Hern{\'{a}}n, Matthew Miller, James~M. Robins, and
  Thomas~S. Richardson.
\newblock Partial identification of the average treatment effect using
  instrumental variables: Review of methods for binary instruments treatments,
  and outcomes.
\newblock \emph{Journal of the American Statistical Association}, 113\penalty0
  (522):\penalty0 933--947, apr 2018.
\newblock \doi{10.1080/01621459.2018.1434530}.
\newblock URL \url{https://doi.org/10.1080%2F01621459.2018.1434530}.

\bibitem[VanderWeele(2015)]{vanderweele2015explanation}
Tyler~J VanderWeele.
\newblock \emph{Explanation in causal inference: methods for mediation and
  interaction}.
\newblock Oxford University Press, 2015.

\bibitem[VanderWeele and Richardson(2012)]{vanderweele2012general}
Tyler~J VanderWeele and Thomas~S Richardson.
\newblock General theory for interactions in sufficient cause models with
  dichotomous exposures.
\newblock \emph{The Annals of Statistics}, 40\penalty0 (4):\penalty0
  2128--2161, 2012.

\bibitem[VanderWeele and Robins(2008)]{vanderweele2008empirical}
Tyler~J VanderWeele and James~M Robins.
\newblock Empirical and counterfactual conditions for sufficient cause
  interactions.
\newblock \emph{Biometrika}, 95\penalty0 (1):\penalty0 49--61, 2008.

\bibitem[VanderWeele et~al.(2012)VanderWeele, Mumford, and
  Schisterman]{vanderweele2012conditioning}
Tyler~J VanderWeele, Sunni~L Mumford, and Enrique~F Schisterman.
\newblock Conditioning on intermediates in perinatal epidemiology.
\newblock \emph{Epidemiology (Cambridge, Mass.)}, 23\penalty0 (1):\penalty0 1,
  2012.

\bibitem[Yerushalmy(1964)]{yerushalmy1964mother}
J~Yerushalmy.
\newblock Mother's cigarette smoking and survival of infant.
\newblock \emph{American Journal of Obstetrics \& Gynecology}, 88\penalty0
  (4):\penalty0 505--518, 1964.

\bibitem[Yerushalmy(2014)]{yerushalmy2014relationship}
J~Yerushalmy.
\newblock The relationship of parents' cigarette smoking to outcome of
  pregnancy--implications as to the problem of inferring causation from
  observed associations.
\newblock \emph{International journal of epidemiology}, 43\penalty0
  (5):\penalty0 1355--1366, 2014.

\bibitem[Zhang and Rubin(2003)]{zhang2003estimation}
Junni~L Zhang and Donald~B Rubin.
\newblock Estimation of causal effects via principal stratification when some
  outcomes are truncated by ``death''.
\newblock \emph{Journal of Educational and Behavioral Statistics}, 28\penalty0
  (4):\penalty0 353--368, 2003.

\end{thebibliography}

\par\null

\section{\texorpdfstring{{Appendix}}{Appendix}}

{\label{987533}}

\subsection{Proofs of Theorems}

{\label{456044}}

\begin{proof}[Proof of Theorem 1]
We prove the contrapositive. Assume that for some \(y\in\left\{0,1\right\}\) and \(m\in\left\{0,1\right\}\) no individual $\omega$ of response type $Y_{1}(\omega)=y,$ $Y_{0}(\omega)=1-y,$ $M_{1}(\omega)=m$ and $M_{0}(\omega)=m$ exists in our population. Then for all individuals $\omega$ in our population $\Omega,$ $I\left\{(Y_{1}(\omega)=y,M_{1}(\omega)=m\right\}+I\left\{Y_{0}(\omega)=1-y,M_{0}(\omega)=m\right\}\leq1.$ Here, $I(\cdot)$ denotes the usual indicator function.

Taking expectations, 
\begin{align*}
E\left[I\left\{Y_{1}(\omega)=y,M_{1}(\omega)=m\right\}+I\left\{Y_{0}(\omega)=1-y,M_{0}(\omega)=m\right\}\right] & \leq & 1\iff \nonumber\\
P\left\{Y_{1}(\omega)=y,M_{1}(\omega)=m\right\}+P\left\{Y_{0}(\omega)=1-y,M_{0}(\omega)=m\right\} & \leq & 1\iff \nonumber\\
P\left\{Y_{1}(\omega)=y,M_{1}(\omega)=m\mid X=1\right\}+P\left\{Y_{0}(\omega)=1-y,M_{0}(\omega)=m\mid X=0\right\} & \leq & 1\iff \nonumber\\
P(Y=y,M=m\mid X=1)+P(Y=1-y,M=m\mid X=0) & \leq & 1.\nonumber
\end{align*}
The second to third line uses $(M_{x},Y_{x})\amalg X$ and the third
to last line uses consistency of counterfactuals. 
\end{proof}

\begin{proof}[Proof of Corollary 2]

From Theorem 1, we have that there exists a non-empty subpopulation $\Omega_{s}\subseteq\Omega$ such that for every $\omega\in\Omega_{s},$ $Y_{1}(\omega)=1,$ $Y_{0}(\omega)=0,$ $M_{1}(\omega)=0$ and $M_{0}(\omega)=0$. From the composition axiom and $M_{1}(\omega)=M_{0}(\omega)=0$, we have that $Y_{1}(\omega)=Y_{1,M_1}(\omega)=Y_{1,M_0}(\omega)$ and $Y_{0}(\omega)=Y_{0,M_0}(\omega)=Y_{0,M_1}(\omega)$, which provides the following results for all $\omega\in\Omega_s$: (1) $Y_{1}(\omega)-Y_{0}(\omega)=1;$ (2) $Y_{1,M_{0}}(\omega)-Y_{0,M_{0}}(\omega)=Y_{1,M_{1}}(\omega)-Y_{0,M_{0}}(\omega)=1;$ (3) $Y_{1,M_{1}}(\omega)-Y_{1,M_{0}}(\omega)=Y_{1,M_{1}}(\omega)-Y_{1,M_{1}}(\omega)=0;$ (4) $Y_{1,M_{1}}(\omega)-Y_{0,M_{1}}(\omega)=Y_{1,M_{1}}(\omega)-Y_{0,M_{0}}(\omega)=1;$ (5) $Y_{0,M_{1}}(\omega)-Y_{0,M_{0}}(\omega)=Y_{0,M_{0}}(\omega)-Y_{0,M_{0}}(\omega)=0;$ (6) $Y_{1,0}(\omega)-Y_{0,0}(\omega)=Y_{1,M_{1}}(\omega)-Y_{0,M_{0}}(\omega)=1;$ (7) $Y_{1}(\omega)-Y_{0}(\omega)=1\text{ and }M_{1}(\omega)=M_{0}(\omega)=0.$
\end{proof}

\begin{proof}[Proof of Proposition 3 and Corollary 4]

Taking expectation of $Y_{1}M_{1}+(1-Y_{0})M_{0}-1,$ $(1-Y_{1})M_{1}+Y_{0}M_{0}-1,$ $Y_{1}(1-M_{1})+(1-Y_{0})(1-M_{0})-1,$ and $(1-Y_{1})(1-M_{1})+Y_{0}(1-M_{0})-1$ gives the result of Proposition 3. A table with the relevant frequencies of individuals $\omega\in\Omega$  is provided in the online supplement. A proof of these results using stochastic counterfactuals is also provided in the online supplement. 

The proof that $$\frac{P \left (Y=y,M=m\mid X=1\right )+P(Y=1-y,M=m\mid X=0)-1}{P \left (M=m\mid X=1\right )-P(M=m-1\mid X=0)}\nonumber$$ is a valid lower bound of $\left(P_c(y,1-y,m,m)-P_c(1-y,y,m,m)\right)/P(M_1=M_0=m)$ is too lengthy for the main text. We provide it in the online supplement. 

We provide a sketch proof regarding the lower bound of $\left(P_c(y,1-y,m,m)-P_c(1-y,y,m,m)\right)/P(M_1=M_0=m)$. Without loss of generality, consider the case where $y=1$ and $m=1$. Note, when $P \left (Y=1,M=1\mid X=1\right )+P(Y=0,M=1\mid X=0)>1,$ then the expression $\left(P_c(1,0,1,1)-P_c(0,1,1,1)\right)/P(M_1=M_0=1)$ can always be  reformulated, using arguments of consistency and randomization, to equal $$A\cdot \frac{P \left (Y=1,M=1\mid X=1\right )+P(Y=0,M=1\mid X=0)-1}{P \left (M=1\mid X=1\right )-P(M=0\mid X=0)},\nonumber$$ where $A$ is a positive constant always greater than equal to $1$. Since constant $A$ is always greater than equal to $1,$ $\left(P_c(1,0,1,1)-P_c(0,1,1,1)\right)/P(M_1=M_0=1)$ is bounded below from $$\frac{P \left (Y=1,M=1\mid X=1\right )+P(Y=0,M=1\mid X=0)-1}{P \left (M=1\mid X=1\right )-P(M=0\mid X=0)}.\nonumber$$

Now, for the last portion of Corollary 4, note $\left(P_c(1,0,1,1)-P_c(0,1,1,1)\right)/P(M_1=M_0=1)$ is equal to $P(Y_1=1,Y_0=0\mid M_1=1,M_0=1)-P(Y_1=0,Y_0=1\mid M_1=1,M_0=1).$ This completes our proof.
\end{proof}

\begin{proof}[Proof of Theorem 5]
We prove the contrapositive. Assume that no individual $\omega$ of response type $Y_{1}(\omega)=1,$ $Y_{0}(\omega)=0,$ $M_{1}(\omega)=0$ and $M_{0}(\omega)=0$ exists in our population. Then for all individuals $\omega$ in our population $\Omega,$\\ $I\left\{Y_{1}(\omega)=1,M_{1}(\omega)=0\right\}-I\left\{Y_{0}(\omega)=1,M_{0}(\omega)=0\right\}\leq 0.$ This last assertion is true after examining a counterfactual table that is provided in the online supplement. Taking expectations, we have 
\begin{align*}
P\left\{Y_{1}(\omega)=1,M_{1}(\omega)=0\right\}-P\left\{Y_{0}(\omega)=1,M_{0}(\omega)=0\right\}& \leq & 0\iff \nonumber\\
P\left\{Y_{1}(\omega)=1,M_{1}(\omega)=0\mid X=1)\right\}-P\left\{Y_{0}(\omega)=1,M_{0}(\omega)=0\mid X=0\right\} & \leq & 0\iff \nonumber \\
P(Y=1,M=0\mid X=1)-P(Y=1,M=0\mid X=0) & \leq & 0.\nonumber
\end{align*}
The first to second line uses $(M_{x},Y_{x})\amalg X$ and the second
to third line uses consistency of counterfactuals. 
\end{proof}

\begin{proof}[Proof of Corollary 6]
Similar to Corollary 2.
\end{proof}

\begin{proof}[Proof of Proposition 7 and Corollary 8]
Taking expectation of $Y_{1}(1-M_{1})-Y_{0}(1-M_{0})$, $Y_{0}M_{0}-Y_{1}M_{1},$ $(1-Y_{1})M_{1}-(1-Y_{0})M_{0}$, and $(1-Y_{1})(1-M_{1})-(1-Y_{0})(1-M_{0})$ gives the result. A table with the relevant frequencies of individuals $\omega\in\Omega$
is provided in the online supplement.

The proof of Corollary 8 is provided next. Under monotonicity, $P(Y=y,M=m\mid X=1-m)-P(Y=y,M=m\mid X=m)$ is lower bound of $P_c(m\{1-y\}+y\{1-m\},ym+\{1-m\}\{1-y\},m,m)-P_c(ym+\{1-m\}\{1-y\},m\{1-y\}+y\{1-m\},m,m),$ and therefore $P(M_1=M_0=m)>0.$ Now, $P(M=m\mid X=1)=P(M_1=m)$ and $P(M=m\mid X=0)=P(M_0=m)$ are both upper bounds to $P(M_1=m,M_0=m).$ We have $P(M=m\mid X=0)\leq P(M=m\mid X=1)$ by the positive monotonicity assumption. Now, dividing a nonnegative lower bound of a numerator by a nonegative upper bound of the denominator will give a lower bound in of the ratio. Choose $P(M=m\mid X=0)$ as the denominator instead of $P(M=m\mid X=1)$ to get greater lower bound. We complete our proof.
\end{proof}

\begin{proof}[Proof of Theorem 9]
For $y\in \{0,1\}$ and $m\in \{0,1\}$, taking expectation of $I(Y_{m}=y)I(M_{m}=m)-I(Y_{1-m}=y)I(M_{1-m}=m)$ and then subtracting $r$ as defined in Theorem 9, gives the result. For the population average principal stratum direct effect, note that $P(M=m\mid X=1)=P(M_1=m)=\sum_{j\in\{0,1\}}P(M_1=m,M_{0}=j).$ Subracting $q$ in the denominator leaves us with $P(M_1=m,M_0=m),$ and hence we have standardized the ratio appropriately.
\end{proof}

\begin{proof}[Proof of Theorem 10]
We prove the contrapositive. Assume that no individual $\omega$ of
response type $Y_{1}(\omega)=1,$ $Y_{0}(\omega)=0,$
$Z_{1}(\omega)=1$ and $Z_{0}(\omega)=0$ exists in our population.
Then for all individuals $\omega$ in our population $\Omega,$\\ $I\left\{Y_{1}(\omega)=y,Z_{1}(\omega)=1\right\}+I\left\{Y_{0}(\omega)=0,Z_{0}(\omega)=0\right\}\leq1.$
Taking expectations, 
\begin{align*}
P\left\{Y_{1}(\omega)=y,Z_{1}(\omega)=1\right\}+P\left\{Y_{0}(\omega)=0,Z_{0}(\omega)=0\right\} & \leq & 1\iff \nonumber \\
P\left\{Y_{1}(\omega)=1,Z_{1}(\omega)=1\mid X=1\right\}+P\left\{Y_{0}(\omega)=0,Z_{0}(\omega)=0\mid  X=0\right\} & \leq & 1\iff \nonumber \\
P(Y=1,Z=1\mid X=1)+P(Y=0,Z=0\mid X=0) & \leq & 1.\nonumber
\end{align*}
The first to second line uses $(Z_{x},Y_{x})\amalg X$ and the second to third line uses consistency of counterfactuals.

To prove the bounding of $P(Y_1=1,Y_0=0,Z_1=1,Z_0=0),$ use the randomization and consistency assumptions. Note,

\begin{align*}
 &P(Y=1,Z=1\mid X=1) + P(Y=0,Z=0\mid X=0)-1 \nonumber\\
 &=P(Y_1=1,Z_1=1)+P(Y_0=0,Z_0=0)-1\nonumber\\
 &=\sum_{(y_0,z_0)\in \{0,1\}^2}P(Y_1=1,Z_1=1,Y_0=y_0,Z_0=z_0)-P((Y_0,Z_0)\neq (0,0))\\
 &=\sum_{(y_0,z_0)\in \{0,1\}^2}P(Y_1=1,Z_1=1,Y_0=y_0,Z_0=z_0)-\sum_{(y_1,z_1)\in \{0,1\}^2}P(Y_1=y_1,Z_1=z_1,(Y_0,Z_0)\neq (0,0))\nonumber\\
 &=P(Y_1=1,Y_0=0,Z_1=1,Z_0=0)-\sum_{(y_1,z_1)\in \{0,1\}^2\setminus \{1,1\}}P(Y_1=y_1,Z_1=z_1,(Y_0,Z_0)\neq (0,0)).\nonumber
\end{align*}
Consequently, $P(Y=1,Z=1\mid X=1)+P(Y=0,Z=0\mid X=0)-1$ is a lower bound of $P(Y_1=1,Y_0=0,Z_1=1,Z_0=0).$
\end{proof}

\begin{proof}[Proof of Proposition 11]
For some $m_a,$ $m_b$ that are subsets of $\mathbf{R}^{n_m}$ and $y_a,$ $y_b$ that are subsets of $\mathbf{R}^{n_y}$, begin with the empirical condition 
\begin{align}
 & P\left(Y\in y_{a},M\in m_{a}\mid X=1\right)+P(Y\not\in y_{b},M\in m_{b}\mid X=0)-1 \label{empab}\\
 & =P\left(Y\in y_{a},M\in m_{a}\mid X=1\right)-P\left(Y\in y_{b}\text{ or }M\not\in m_{b}\mid X=0\right)\nonumber\\
 & =P\left(Y\in y_{a},M\in m_{a}\mid X=1\right)-P\left(Y\in y_{b},M\in m_{b}\mid X=0\right)\nonumber\\
 &\hspace{1 em} -P\left(Y\not\in y_{b},M\not\in m_{b}\mid X=0\right)-P\left(Y\in y_{b},M\not\in m_{b}\mid X=0\right)\nonumber
\end{align}
 Under randomization and consistency of counterfactuals, the empirical
condition (\ref{empab}) is equal to 
\begin{align*}
 & P\left(Y_{1}\in y_{a},M_{1}\in m_{a}\right)-P\left(Y_{0}\in y_{b},M_{0}\in m_{b}\right)\nonumber-P\left(Y_{0}\not\in y_{b},M_{0}\not\in m_{b}\right)-P\left(Y_{0}\in y_{b},M_{0}\not\in m_{b}\right).
\end{align*}
Use the law of total probability. The empirical
condition (\ref{empab}) is equal to 
\begin{align*}
 & P\left(Y_{1}\in y_{a},Y_{0}\in y_{b},M_{1}\in m_{a},M_{0}\in m_{b}\right)+P\left(Y_{1}\in y_{a},Y_{0}\not\in y_{b},M_{1}\in m_{a},M_{0}\in m_{b}\right)\nonumber\\
 & +P\left(Y_{1}\in y_{a},Y_{0}\in y_{b},M_{1}\in m_{a},M_{0}\not\in m_{b}\right)+P\left(Y_{1}\not\in y_{a},Y_{0}\in y_{b},M_{1}\in m_{a},M_{0}\not\in m_{b}\right)\nonumber\\
 & -\left\{ P\left(Y_{1}\in y_{a},Y_{0}\in y_{b},M_{1}\in m_{a},M_{0}\in m_{b}\right)+P\left(Y_{1}\not\in y_{a},Y_{0}\in y_{b},M_{1}\in m_{a},M_{0}\in m_{b}\right)\right.\nonumber\\
 & +P\left(Y_{1}\in y_{a},Y_{0}\in y_{b},M_{1}\not\in m_{a},M_{0}\in m_{b}\right)+P\left(Y_{1}\not\in y_{a},Y_{0}\in y_{b},M_{1}\not\in m_{a},M_{0}\in m_{b}\right)\nonumber\\
 & +P\left(Y_{1}\in y_{a},Y_{0}\not\in y_{b},M_{1}\in m_{a},M_{0}\not\in m_{b}\right)+P\left(Y_{1}\not\in y_{a},Y_{0}\not\in y_{b},M_{1}\in m_{a},M_{0}\not\in m_{b}\right)\nonumber\\
 & +P\left(Y_{1}\in y_{a},Y_{0}\not\in y_{b},M_{1}\not\in m_{a},M_{0}\not\in m_{b}\right)+P\left(Y_{1}\not\in y_{a},Y_{0}\not\in y_{b},M_{1}\not\in m_{a},M_{0}\not\in m_{b}\right)\nonumber\\
 & +P\left(Y_{1}\in y_{a},Y_{0}\in y_{b},M_{1}\in m_{a},M_{0}\not\in m_{b}\right)+P\left(Y_{1}\not\in y_{a},Y_{0}\in y_{b},M_{1}\in m_{a},M_{0}\not\in m_{b}\right)\nonumber\\
 & \left.+P\left(Y_{1}\in y_{a},Y_{0}\in y_{b},M_{1}\not\in m_{a},M_{0}\not\in m_{b}\right)+P\left(Y_{1}\not\in y_{a},Y_{0}\in y_{b},M_{1}\not\in m_{a},M_{0}\not\in m_{b}\right)\right\} .
\end{align*}
 Simplifying this last expression, we have that the empirical condition (\ref{empab})
is equal to 
\begin{align}
 & P(Y_{1}\in y_{a},Y_{0}\not\in y_{b},M_{1}\in m_{a},M_{0}\in m_{b})-P\left(Y_{1}\not\in y_{a},Y_{0}\in y_{b},M_{1}\in m_{a},M_{0}\in m_{b}\right)\nonumber\\
 & -\left\{ P\left(Y_{1}\in y_{a},Y_{0}\in y_{b},M_{1}\not\in m_{a},M_{0}\in m_{b}\right)+P\left(Y_{1}\not\in y_{a},Y_{0}\in y_{b},M_{1}\not\in m_{a},M_{0}\in m_{b}\right)\right.\nonumber\\
 & +P\left(Y_{1}\in y_{a},Y_{0}\not\in y_{b},M_{1}\in m_{a},M_{0}\not\in m_{b}\right)+P\left(Y_{1}\not\in y_{a},Y_{0}\not\in y_{b},M_{1}\in m_{a},M_{0}\not\in m_{b}\right)\\
 & +P\left(Y_{1}\in y_{a},Y_{0}\not\in y_{b},M_{1}\not\in m_{a},M_{0}\not\in m_{b}\right)+P\left(Y_{1}\not\in y_{a},Y_{0}\not\in y_{b},M_{1}\not\in m_{a},M_{0}\not\in m_{b}\right)\nonumber\\
 & \left.+P\left(Y_{1}\in y_{a},Y_{0}\in y_{b},M_{1}\not\in m_{a},M_{0}\not\in m_{b}\right)+P\left(Y_{1}\not\in y_{a},Y_{0}\in y_{b},M_{1}\not\in m_{a},M_{0}\not\in m_{b}\right).\right\} \nonumber
\end{align}
 Given that probability is greater than equal to zero and less than
equal to 1, the empirical condition (\ref{empab}) is a lower bound on $P(Y_{1}\in y_{a},Y_{0}\not\in y_{b},M_{1}\in m_{a},M_{0}\in m_{b})-P\left(Y_{1}\not\in y_{a},Y_{0}\in y_{b},M_{1}\in m_{a},M_{0}\in m_{b}\right).$ 

If we take $y_{b}=y_{a}$ and $m_{b}=m_{a},$ then $P(Y\in y_{a},M\in m_{a}\mid X=1)+P(Y\not\in y_{a},M\in m_{a}\mid X=0)-1$ serves as a lower bound on $P(Y_{1}\in y_{a},Y_{0}\not\in y_{a},M_{1}\in m_{a},M_{0}\in m_{a})-P\left(Y_{1}\not\in y_{a},Y_{0}\in y_{a},M_{1}\in m_{a},M_{0}\in m_{a}\right).$ 

Additionally, 
\begin{align*}
 &P (Y\in y_{a},M\in m_{a}\mid X=1)+P(Y\not\in y_{a},M\in m_{a}\mid X=0)-1\\
 & = P(Y_{1}\in y_{a},Y_{0}\not\in y_{a},M_{1}\in m_{a},M_{0}\in m_{a})-P\left(Y_{1}\not\in y_{a},Y_{0}\in y_{a},M_{1}\in m_{a},M_{0}\in m_{a}\right)\nonumber\\
 & \hspace{1 em} -\left\{ P\left(Y_{1}\in y_{a},Y_{0}\in y_{a},M_{1}\not\in m_{a},M_{0}\in m_{a}\right)+P\left(Y_{1}\not\in y_{a},Y_{0}\in y_{a},M_{1}\not\in m_{a},M_{0}\in m_{a}\right)\right.\nonumber\\
 & \hspace{1 em} +P\left(Y_{1}\in y_{a},Y_{0}\not\in y_{a},M_{1}\in m_{a},M_{0}\not\in m_{a}\right)+P\left(Y_{1}\not\in y_{a},Y_{0}\not\in y_{a},M_{1}\in m_{a},M_{0}\not\in m_{a}\right)\nonumber\\
 & \hspace{1 em} +P\left(Y_{1}\in y_{a},Y_{0}\not\in y_{a},M_{1}\not\in m_{a},M_{0}\not\in m_{a}\right)+P\left(Y_{1}\not\in y_{a},Y_{0}\not\in y_{a},M_{1}\not\in m_{a},M_{0}\not\in m_{a}\right)\nonumber\\
 & \hspace{1 em} \left.+P\left(Y_{1}\in y_{a},Y_{0}\in y_{a},M_{1}\not\in m_{a},M_{0}\not\in m_{a}\right)+P\left(Y_{1}\not\in y_{a},Y_{0}\in y_{a},M_{1}\not\in m_{a},M_{0}\not\in m_{a}\right).\right\}
\end{align*}
This is true if and only if 
\begin{align*}
 & \frac{P(Y_{1}\in y_{a},Y_{0}\not\in y_{a},M_{1}\in m_{a},M_{0}\in m_{a})-P\left(Y_{1}\not\in y_{a},Y_{0}\in y_{a},M_{1}\in m_{a},M_{0}\in m_{a}\right)}{P(M_{1}\in m_{a},M_{0}\in m_{a})}\nonumber\\
 & =\frac{P(Y\in y_{a},M\in m_{a}\mid X=1)+P(Y\not\in y_{a},M\in m_{a}\mid X=0)-1}{P(M_{1}\in m_{a},M_{0}\in m_{a})}\nonumber\\
 & +\frac{P\left(Y_{1}\in y_{a},Y_{0}\in y_{a},M_{1}\not\in m_{a},M_{0}\in m_{a}\right)+P\left(Y_{1}\not\in y_{a},Y_{0}\in y_{a},M_{1}\not\in m_{a},M_{0}\in m_{a}\right)}{P(M_{1}\in m_{a},M_{0}\in m_{a})}\nonumber\\
 & +\frac{P\left(Y_{1}\in y_{a},Y_{0}\not\in y_{a},M_{1}\in m_{a},M_{0}\not\in m_{a}\right)+P\left(Y_{1}\not\in y_{a},Y_{0}\not\in y_{a},M_{1}\in m_{a},M_{0}\not\in m_{a}\right)}{P(M_{1}\in m_{a},M_{0}\in m_{a})}\nonumber\\
 & +\frac{P\left(Y_{1}\in y_{a},Y_{0}\not\in y_{a},M_{1}\not\in m_{a},M_{0}\not\in m_{a}\right)+P\left(Y_{1}\not\in y_{a},Y_{0}\not\in y_{a},M_{1}\not\in m_{a},M_{0}\not\in m_{a}\right)}{P(M_{1}\in m_{a},M_{0}\in m_{a})}\nonumber\\
 & +\frac{P\left(Y_{1}\in y_{a},Y_{0}\in y_{a},M_{1}\not\in m_{a},M_{0}\not\in m_{a}\right)+P\left(Y_{1}\not\in y_{a},Y_{0}\in y_{a},M_{1}\not\in m_{a},M_{0}\not\in m_{a}\right)}{P(M_{1}\in m_{a},M_{0}\in m_{a})}.
\end{align*}
Using the identity, 
\begin{align}
 & P(M_{1}\not\in m_{a},M_{0}\not\in m_{a})\nonumber\\
 & =P\left(Y_{1}\in y_{a},Y_{0}\not\in y_{a},M_{1}\not\in m_{a},M_{0}\not\in m_{a}\right)+P\left(Y_{1}\not\in y_{a},Y_{0}\not\in y_{a},M_{1}\not\in m_{a},M_{0}\not\in m_{a}\right) \label{notm1m0}\\
 & \hspace{1 em} +P\left(Y_{1}\in y_{a},Y_{0}\in y_{a},M_{1}\not\in m_{a},M_{0}\not\in m_{a}\right)+P\left(Y_{1}\not\in y_{a},Y_{0}\in y_{a},M_{1}\not\in m_{a},M_{0}\not\in m_{a}\right),\nonumber
\end{align}
this equation for the standardized principal stratum direct effect simplifies
to
\begin{align*}
 & \frac{P(Y_{1}\in y_{a},Y_{0}\not\in y_{a},M_{1}\in m_{a},M_{0}\in m_{a})-P\left(Y_{1}\not\in y_{a},Y_{0}\in y_{a},M_{1}\in m_{a},M_{0}\in m_{a}\right)}{P(M_{1}\in m_{a},M_{0}\in m_{a})}\nonumber\\
 & =\frac{P(Y\in y_{a},M\in m_{a}\mid X=1)+P(Y\not\in y_{a},M\in m_{a}\mid X=0)-1}{P(M_{1}\in m_{a},M_{0}\in m_{a})}\nonumber\\
 & +\frac{P\left(Y_{1}\in y_{a},Y_{0}\in y_{a},M_{1}\not\in m_{a},M_{0}\in m_{a}\right)+P\left(Y_{1}\not\in y_{a},Y_{0}\in y_{a},M_{1}\not\in m_{a},M_{0}\in m_{a}\right)}{P(M_{1}\in m_{a},M_{0}\in m_{a})}\nonumber\\
 & +\frac{P\left(Y_{1}\in y_{a},Y_{0}\not\in y_{a},M_{1}\in m_{a},M_{0}\not\in m_{a}\right)+P\left(Y_{1}\not\in y_{a},Y_{0}\not\in y_{a},M_{1}\in m_{a},M_{0}\not\in m_{a}\right)}{P(M_{1}\in m_{a},M_{0}\in m_{a})}\nonumber\\
 & +\frac{P(M_{1}\not\in m_{a},M_{0}\not\in m_{a})}{P(M_{1}\in m_{a},M_{0}\in m_{a})}.\nonumber
\end{align*}
Denote 
\begin{align*}
x & =P\left(Y_{1}\not\in y_{a},Y_{0}\not\in y_{a},M_{1}\in m_{a},M_{0}\not\in m_{a}\right)+P\left(Y_{1}\in y_{a},Y_{0}\not\in y_{a},M_{1}\in m_{a},M_{0}\not\in m_{a}\right)\nonumber\\
 &\hspace{1 em} +P\left(Y_{1}\not\in y_{a},Y_{0}\in y_{a},M_{1}\not\in m_{a},M_{0}\in m_{a}\right)+P\left(Y_{1}\in y_{a},Y_{0}\in y_{a},M_{1}\not\in m_{a},M_{0}\in m_{a}\right)\nonumber
\end{align*}
 out of space considerations. Note, $x\in[0,1).$ With some algebra, the standardized principal stratum direct effect is equal to
\begin{align*}
 & \frac{P(Y_{1}\in y_{a},Y_{0}\not\in y_{a},M_{1}\in m_{a},M_{0}\in m_{a})-P\left(Y_{1}\not\in y_{a},Y_{0}\in y_{a},M_{1}\in m_{a},M_{0}\in m_{a}\right)}{P(M_{1}\in m_{a},M_{0}\in m_{a})}\\
 & =\frac{P(Y\in y_{a},M\in m_{a}\mid X=1)+P(Y\not\in y_{a},M\in m_{a}\mid X=0)-1}{P(M_{1}\in m_{a},M_{0}\in m_{a})}\nonumber\\
 & +\frac{x+P(M_{1}\not\in m_{a},M_{0}\not\in m_{a})}{P(M_{1}\in m_{a},M_{0}\in m_{a})}.\nonumber\\
 & =\frac{P(Y\in y_{a},M\in m_{a}\mid X=1)+P(Y\not\in y_{a},M\in m_{a}\mid X=0)-1}{P(M_{1}\in m_{a},M_{0}\in m_{a})-P(M_{1}\not\in m_{a},M_{0}\not\in m_{a})}\\
 &\hspace{1.5 em}\times \left(\frac{P(M_{1}\in m_{a},M_{0}\in m_{a})-P(M_{1}\not\in m_{a},M_{0}\not\in m_{a})}{P(M_{1}\in m_{a},M_{0}\in m_{a})}\right.\nonumber\\
 &\hspace{1.6 em} \left.+\frac{\left(x+P(M_{1}\not\in m_{a},M_{0}\not\in m_{a})\right)}{P(M_{1}\in m_{a},M_{0}\in m_{a})}\frac{\left(P(M_{1}\in m_{a},M_{0}\in m_{a})-P(M_{1}\not\in m_{a},M_{0}\not\in m_{a})\right)}{P(Y\in y_{a},M\in m_{a}\mid X=1)+P(Y\not\in y_{a},M\in m_{a}\mid X=0)-1}\right).\nonumber
\end{align*}
We want to show that the next expression
\begin{align}
\left(\frac{P(M_{1}\in m_{a},M_{0}\in m_{a})}{P(M_{1}\in m_{a},M_{0}\in m_{a})-P(M_{1}\not\in m_{a},M_{0}\not\in m_{a})}\right.\nonumber\\
\left.+\frac{\left(x+P(M_{1}\not\in m_{a},M_{0}\not\in m_{a})\right)}{P(M_{1}\in m_{a},M_{0}\in m_{a})}\frac{\left(P(M_{1}\in m_{a},M_{0}\in m_{a})-P(M_{1}\not\in m_{a},M_{0}\not\in m_{a})\right)}{P(Y\in y_{a},M\in m_{a}\mid X=1)+P(Y\not\in y_{a},M\in m_{a}\mid X=0)-1}\right).\ \label{expA}
\end{align}
is greater than or equal to one to show that $\frac{P(Y\in y_{a},M\in m_{a}\mid X=1)+P(Y\not\in y_{a},M\in m_{a}\mid X=0)-1}{P(M_{1}\in m_{a},M_{0}\in m_{a})-P(M_{1}\not\in m_{a},M_{0}\not\in m_{a})}$
is a lower bound of 
\[
\frac{P(Y_{1}\in y_{a},Y_{0}\not\in y_{a},M_{1}\in m_{a},M_{0}\in m_{a})-P\left(Y_{1}\not\in y_{a},Y_{0}\in y_{a},M_{1}\in m_{a},M_{0}\in m_{a}\right)}{P(M_{1}\in m_{a},M_{0}\in m_{a})}.
\nonumber\]
Through inspection, expression (\ref{expA}) is greater than or equal to one if
and only if 
\[
\frac{\left(P(M_{1}\in m_{a},M_{0}\in m_{a})-P(M_{1}\not\in m_{a},M_{0}\not\in m_{a})\right)}{P(Y\in y_{a},M\in m_{a}\mid X=1)+P(Y\not\in y_{a},M\in m_{a}\mid X=0)-1}\geq1.\nonumber
\]
Identity (\ref{notm1m0}) provides the equality
\begin{align*}
 & P(Y\in y_{a},M\in m_{a}\mid X=1)+P(Y\not\in y_{a},M\in m_{a}\mid X=0)-1\nonumber\\
 & =P(Y_{1}\in y_{a},Y_{0}\not\in y_{a},M_{1}\in m_{a},M_{0}\in m_{a})-P\left(Y_{1}\not\in y_{a},Y_{0}\in y_{a},M_{1}\in m_{a},M_{0}\in m_{a}\right)\nonumber\\
 &\hspace{1 em} -\left\{ P\left(Y_{1}\in y_{a},Y_{0}\in y_{a},M_{1}\not\in m_{a},M_{0}\in m_{a}\right)+P\left(Y_{1}\not\in y_{a},Y_{0}\in y_{a},M_{1}\not\in m_{a},M_{0}\in m_{a}\right)\right.\nonumber\\
 &\hspace{1 em} +P\left(Y_{1}\in y_{a},Y_{0}\not\in y_{a},M_{1}\in m_{a},M_{0}\not\in m_{a}\right)+P\left(Y_{1}\not\in y_{a},Y_{0}\not\in y_{a},M_{1}\in m_{a},M_{0}\not\in m_{a}\right)\nonumber\\
 &\hspace{1 em} \left.+P\left(M_{1}\not\in m_{a},M_{0}\not\in m_{a}\right)\right\} .
\end{align*}
The fraction
\[\nonumber
\frac{\left(P(M_{1}\in m_{a},M_{0}\in m_{a})-P(M_{1}\not\in m_{a},M_{0}\not\in m_{a})\right)}{P(Y\in y_{a},M\in m_{a}\mid X=1)+P(Y\not\in y_{a},M\in m_{a}\mid X=0)-1}\nonumber
\]
is greater than equal to one, if and only if the numerator is greater
than the denominator. Here, we have that the denominator is positive
as a condition of our Theorem. The only positive term in the denominator
is $P(Y_{1}\in y_{a},Y_{0}\not\in y_{a},M_{1}\in m_{a},M_{0}\in m_{a}),$
which is less than equal to $P(M_{1}\in m_{a},M_{0}\in m_{a}).$ Also,
the numerator is greater than zero, as 
\begin{align*}\nonumber
P(M_{1}\in m_{a},M_{0}\in m_{a})-P(M_{1}\not\in m_{a},M_{0}\not\in m_{a})\geq P(Y_{1}\in y_{a},Y_{0}\not\in y_{a},M_{1}\in m_{a},M_{0}\in m_{a})-P\left(M_{1}\not\in m_{a},M_{0}\not\in m_{a}\right).
\end{align*}
Therefore, 
\[\nonumber
\frac{\left(P(M_{1}\in m_{a},M_{0}\in m_{a})-P(M_{1}\not\in m_{a},M_{0}\not\in m_{a})\right)}{P(Y\in y_{a},M\in m_{a}\mid X=1)+P(Y\not\in y_{a},M\in m_{a}\mid X=0)-1}\geq1.\nonumber
\]
Consequently, our expression (\ref{expA}) is greater than or equal to one, which in turn implies that 
\begin{align*}
 & \frac{P(Y\in y_{a},M\in m_{a}\mid X=1)+P(Y\not\in y_{a},M\in m_{a}\mid X=0)-1}{P(M_{1}\in m_{a},M_{0}\in m_{a})-P(M_{1}\not\in m_{a},M_{0}\not\in m_{a})}\nonumber\\
 & \leq\frac{P(Y_{1}\in y_{a},Y_{0}\not\in y_{a},M_{1}\in m_{a},M_{0}\in m_{a})-P\left(Y_{1}\not\in y_{a},Y_{0}\in y_{a},M_{1}\in m_{a},M_{0}\in m_{a}\right)}{P(M_{1}\in m_{a},M_{0}\in m_{a})}\nonumber
\end{align*}
 Therefore, 
\begin{align*}
\frac{P(Y\in y_{a},M\in m_{a}\mid X=1)+P(Y\not\in y_{a},M\in m_{a}\mid X=0)-1}{P(M_{1}\in m_{a},M_{0}\in m_{a})-P(M_{1}\not\in m_{a},M_{0}\not\in m_{a})}\nonumber
\end{align*}
 is a lower bound on
\begin{align*}
\frac{P(Y_{1}\in y_{a},Y_{0}\not\in y_{a},M_{1}\in m_{a},M_{0}\in m_{a})-P\left(Y_{1}\not\in y_{a},Y_{0}\in y_{a},M_{1}\in m_{a},M_{0}\in m_{a}\right)}{P(M_{1}\in m_{a},M_{0}\in m_{a})-P(M_{1}\not\in m_{a},M_{0}\not\in m_{a})},\nonumber
\end{align*}
when $P(Y\in y_{a},M\in m_{a}\mid X=1)+P(Y\not\in y_{a},M\in m_{a}\mid X=0)-1>0.$ 

Now, note 
\begin{align*}
 & P(M_{1}\in m_{a},M_{0}\in m_{a})-P(M_{1}\not\in m_{a},M_{0}\not\in m_{a})\nonumber\\
 & =P(M_{1}\in m_{a},M_{0}\in m_{a})+P(M_{1}\in m_{a},M_{0}\not\in m_{a})-P(M_{1}\in m_{a},M_{0}\not\in m_{a})-P(M_{1}\not\in m_{a},M_{0}\not\in m_{a})\\
 & =P(M_{1}\in m_{a})-P(M_{0}\not\in m_{a})\nonumber\\
 & =P(M\in m_{a}\mid X=1)-P(M\not\in m_{a}\mid X=0).
\end{align*}
 Consequently, 
\[
\frac{P(Y\in y_{a},M\in m_{a}\mid X=1)+P(Y\not\in y_{a},M\in m_{a}\mid X=0)-1}{P(M\in m_{a}\mid X=1)-P(M\not\in m_{a}\mid X=0)}\nonumber
\]
is a lower bound of 
\begin{align*}
\frac{P(Y_{1}\in y_{a},Y_{0}\not\in y_{a},M_{1}\in m_{a},M_{0}\in m_{a})-P\left(Y_{1}\not\in y_{a},Y_{0}\in y_{a},M_{1}\in m_{a},M_{0}\in m_{a}\right)}{P(M_{1}\in m_{a},M_{0}\in m_{a})}.\nonumber
\end{align*}
in the case when $P(Y\in y_{a},M\in m_{a}\mid X=1)+P(Y\not\in y_{a},M\in m_{a}\mid X=0)-1>0.$
\end{proof}

\begin{proof}[Proposition 12]
From the randomization assumption, consistency of counterfactuals, and an application of the law of total probability,
\begin{align*}
 & P(Y\not\in y_{a},M\in m_{a}\mid X=0)-P(Y\not\in y_{a},M\in m_{a}\mid X=1)\\
 & =P(Y_{0}\not\in y_{a},M_{0}\in m_{a})-P(Y_{1}\in y_{a},M_{1}\in m_{a})\\
 & =P\left(Y_{0}\not\in y_{a},M_{1}\in m_{a},M_{0}\in m_{a}\right)+P(Y_{0}\not\in y_{a},M_{1}\not\in m_{a},M_{0}\in m_{a})\\
 & \hspace{1 em}-P(Y_{1}\not\in y_{a},M_{1}\in m_{a},M_{0}\in m_{a})-P(Y_{1}\not\in y_{a},M_{1}\in m_{a},M_{0}\not\in m_{a}).
\end{align*}
Another application of the law of total probability, we have 
\begin{align*}
 & P(Y\not\in y_{a},M\in m_{a}\mid X=0)-P(Y\not\in y_{a},M\in m_{a}\mid X=1)\\
 & =P(Y_{1}\in y_{a},Y_{0}\not\in y_{a},M_{1}\in m_{a},M_{0}\in m_{a})+P(Y_{1}\not\in y_{a},Y_{0}\not\in y_{a},M_{1}\in m_{a},M_{0}\in m_{a})\\
 &\hspace{1 em} +P(Y_{1}\in y_{a},Y_{0}\not\in y_{a},M_{1}\not\in m_{a},M_{0}\in m_{a}) +P(Y_{1}\not\in y_{a},Y_{0}\not\in y_{a},M_{1}\not\in m_{a},M_{0}\in m_{a})\\
 &\hspace{1 em} -P(Y_{1}\not\in y_{a},Y_{0}\in y_{a},M_{1}\in m_{a},M_{0}\in m_{a}) -P(Y_{1}\not\in y_{a},Y_{0}\not\in y_{a},M_{1}\in m_{a},M_{0}\in m_{a})\\
 &\hspace{1 em} -P(Y_{1}\not\in y_{a},Y_{0}\in y_{a},M_{1}\in m_{a},M_{0}\not\in m_{a}) -P(Y_{1}\not\in y_{a},Y_{0}\not\in y_{a},M_{1}\in m_{a},M_{0}\not\in m_{a}).
\end{align*}
Simplifying this equality, 
\begin{align*}
 & P(Y\notin y_{a},M\in m_{a}\mid X=0)-P(Y\not\in y_{a},M\in m_{a}\mid X=1)\\
 & =P(Y_{1}\in y_{a},Y_{0}\not\in y_{a},M_{1}\in m_{a},M_{0}\in m_{a})-P(Y_{1}\not\in y_{a},Y_{0}\in y_{a},M_{1}\in m_{a},M_{0}\in m_{a})\\
 &\hspace{1 em} +\left\{ P(Y_{1}\not\in y_{a},Y_{0}\not\in y_{a},M_{1}\not\in m_{a},M_{0}\in m_{a})+P(Y_{1}\in y_{a},Y_{a}\not\in y_{a},M_{1}\not\in m_{a},M_{0}\in m_{a})\right.\\
 &\hspace{1 em} \left. -P(Y_{1}\not\in y_{a},Y_{0}\not\in y_{a},M_{1}\in m_{a},M_{0}\not\in m_{a})-P(Y_{1}\not\in y_{a},Y_{0}\in y_{a},M_{1}\in m_{a},M_{0}\not\in m_{a})\right\}.
\end{align*}
The monotonicity assumption that there is no individual
of response type $M_{1}(\omega)\not\in m_{a}$ and $M_{0}(\omega)\in m_{a}$
and $P(Y\in y_{a},M\in m_{a}\mid X=0)>P(Y\in y_{a},M\in m_{a}\mid X=1),$
yields 
\begin{align*}
 & \left\{ P\left(Y_{1}\not\in y_{a},Y_{0}\not\in Y_{a},M_{1}\not\in m_{a},M_{0}\in m_{a}\right)\right.+P(Y_{1}\in y_{a},Y_{0}\not\in y_{a},M_{1}\not\in m_{a},M_{0}\in m_{a})\\
 & -P(Y_{1}\not\in y_{a},Y_{0}\not\in y_{a},M_{1}\in m_{a},M_{0}\not\in m_{a})\left.-P\left(Y_{1}\not\in y_{a},Y_{0}\in y_{a},M_{1}\in m_{a},M_{0}\not\in m_{a}\right)\right\} \leq0,
\end{align*}
 and $P(Y\not\in y_{a},M\in m_{a}\mid X=0)-P(Y\not\in y_{a},M\in m_{a}\mid X=1)$ is less than $P(Y_{1}\in y_{a},Y_{0}\not\in y_{a},M_{1}\in m_{a},M_{0}\in m_{a}) -P(Y\not\in y_{a},Y_{0}\in y_{a},M_{1}\in m_{a},M_{0}\in m_{a}).$ We also have that $P(M\in m_{a}\mid X=0)$ is an upper bound on $P(M_{1}\in m_{a},M_{0}\in m_{a}).$ 
 
 Provided the numerator is positive,  
$$\frac{\left(P(Y\not\in y_{a},M\in m_{a}\mid X=0)-P(Y\not\in y_{a},M\in m_{a}\mid X=1)\right)}{P(M\in m_{a}\mid X=0)}$$
 is a lower bound on 
\[\frac{\left(\begin{array}{c}
P(Y_{1}\in y_{a},Y_{0}\not\in y_{a},M_{1}\in m_{a},M_{0}\in m_{a})\\
-P(Y_{1}\not\in y_{a},Y_{0}\in y_{a},M_{1}\in m_{a},M_{0}\in m_{a})
\end{array}\right)}{P(M_{1}\in m_{a},M_{0}\in m_{a})},\]
and equivalently $P(Y_{1}\in y_{a},Y_{0}\not\in y_{a}\mid M_{1}\in m_{a},M_{0}\in m_{a})-P(Y_{1}\not\in y_{a},Y_{0}\in y_{a}\mid M_{1}\in m_{a},M_{0}\in m_{a}).$
This completes our proof.
\end{proof}

\begin{proof}[Proof of Proposition 13]
From the proof of Proposition 12 above, we had 
\begin{align*}
 & P(Y\notin y_{a},M\in m_{a}\mid X=0)-P(Y\not\in y_{a},M\in m_{a}\mid X=1)\\
 & =P(Y_{1}\in y_{a},Y_{0}\not\in y_{a},M_{1}\in m_{a},M_{0}\in m_{a})-P(Y_{1}\not\in y_{a},Y_{0}\in y_{a},M_{1}\in m_{a},M_{0}\in m_{a})\\
 & \hspace{1 em} +\left\{ P(Y_{1}\not\in y_{a},Y_{0}\not\in y_{a},M_{1}\not\in m_{a},M_{0}\in m_{a})+P(Y_{1}\in y_{a},Y_{a}\not\in y_{a},M_{1}\not\in m_{a},M_{0}\in m_{a})\right.\\
 & \hspace{1 em} \left.-P(Y_{1}\not\in y_{a},Y_{0}\not\in y_{a},M_{1}\in m_{a},M_{0}\not\in m_{a})-P(Y_{1}\not\in y_{a},Y_{0}\in y_{a},M_{1}\in m_{a},M_{0}\not\in m_{a})\right\}.
\end{align*}
Thus, $P(Y\not\in y_{a},M\in m_{a}\mid X=0)-P(Y\not\in y_{a},M\in m_{a}\mid X=1)-r,$ where 
\begin{align*}
r= & \left\{ P(Y_{1}\not\in y_{a},Y_{0}\not\in y_{a},M_{1}\not\in m_{a},M_{0}\in m_{a})\right.\\
 &\hspace{1 em} +P(Y_{1}\in y_{a},Y_{a}\not\in y_{a},M_{1}\not\in m_{a},M_{0}\in m_{a})\\
 &\hspace{1 em} -P(Y_{1}\not\in y_{a},Y_{0}\not\in y_{a},M_{1}\in m_{a},M_{0}\not\in m_{a})\\
 &\hspace{1 em} \left.-P(Y_{1}\not\in y_{a},Y_{0}\in y_{a},M_{1}\in m_{a},M_{0}\not\in m_{a})\right\}.,
\end{align*}
 is equal to $P(Y_{1}\not\in y_{a},Y_{0}\in y_{a},M_{1}\in m_{a},M_{0}\in m_{a})-P(Y\in y_{a},Y_{0}\not\in y_{a},M_{1}\in m_{a},M_{0}\in m_{a}).$
Now, under randomization and consistency of counterfactuals,
$$P(M\in m_{a}\mid X=0)=P(M_{0}\in m_{a})=P(M_{1}\in m_{a},M_{0}\in m_{a})+P(M_{1}\not\in m_{a},M_{0}\in m_{a}).$$
 Therefore, 
\[
P(M_{1}\in m_{a},M_{0}\in m_{a})=P(M\in m_{a}\mid X=0)-P(M_{1}\not\in m_{a},M_{0}\in m_{a}).
\]
 This means that the expression 
\[
\frac{P(Y_{1}\in y_{a},Y_{0}\not\in y_{a},M_{1}\in m_{a},M_{0}\in m_{a})-P(Y_{1}\not\in y_{a},Y_{0}\in y_{a},M_{1}\in m_{a},M_{0}\in m_{a})}{P(M\in m_{a}\mid X=0)-P(M_{1}\not\in m_{a},M_{0}\in m_{a})}
\]
 is equal to 
\begin{align*}
 & \frac{P(Y\not\in y_{a},M\in m_{a}\mid X=0)-P(Y\not\in y_{a},M\in m_{a}\mid X=1)-r}{P(M\in m_{a}\mid X=0)-q}\\
\end{align*}
where $q=P(M_{1}\not\in m_{a},M_{0}\in m_{a})$ . This completes the
proof.
\end{proof}

\end{document}